# A PROBLEM IN ONE-DIMENSIONAL DIFFUSION-LIMITED AGGREGATION (DLA) AND POSITIVE RECURRENCE OF MARKOV CHAINS

By Harry Kesten and Vladas Sidoravicius

*Cornell University and IMPA, Brasil*

We consider the following problem in one-dimensional diffusion-limited aggregation (DLA). At time $t$, we have an "aggregate" consisting of $\mathbb{Z} \cap [0, R(t)]$ [with $R(t)$ a positive integer]. We also have $N(i,t)$ particles at $i, i > R(t)$. All these particles perform independent continuous-time symmetric simple random walks until the first time $t' > t$ at which some particle tries to jump from $R(t)+1$ to $R(t)$. The aggregate is then increased to the integers in $[0, R(t')] = [0, R(t)+1]$ [so that $R(t') = R(t)+1$] and all particles which were at $R(t)+1$ at time $t'-$ are removed from the system. The problem is to determine how fast $R(t)$ grows as a function of $t$ if we start at time 0 with $R(0) = 0$ and the $N(i,0)$ i.i.d. Poisson variables with mean $\mu > 0$. It is shown that if $\mu < 1$, then $R(t)$ is of order $\sqrt{t}$, in a sense which is made precise. It is conjectured that $R(t)$ will grow linearly in $t$ if $\mu$ is large enough.

**1. Introduction.** Before we begin the discussion of the speed at which the aggregate in the diffusion-limited aggregation (DLA) model grows, we explain how we came to this problem from studying another growth model. In Kesten and Sidoravicius (2005), we studied the following model for the spread of an infection. There is a "gas" of particles, each of which performs a continuous-time simple random walk on $\mathbb{Z}^d$ with jump rate $D_A$. These particles are called *A-particles* and move independently of each other. They are regarded as healthy individuals. We assume that we start the system with $N_A(x, 0-)$ $A$-particles at $x$ and that the $N_A(x, 0-)$, $x \in \mathbb{Z}^d$, are i.i.d., mean-$\mu_A$ Poisson random variables. In addition, there are *B-particles* which perform continuous-time simple random walks with jump rate $D_B$. We start with a finite number of $B$-particles in the system at time 0. $B$-particles are interpreted as infected individuals. The $B$-particles move independently of each other. The only interaction is that when a $B$-particle and an $A$-particle coincide, the latter instantaneously turns into a $B$-particle.

In Kesten and Sidoravicius (2005), we investigated how fast the infection spreads. Specifically, if $\widetilde{B}(t) := \{x \in \mathbb{Z}^d : \text{a } B\text{-particle visits } x \text{ during } [0,t]\}$ and $B(t) = \widetilde{B}(t) + [-1/2, 1/2]^d$, then we investigated the asymptotic behavior of $B(t)$. The principal result in Kesten and Sidoravicius (2005) states that if $D_A = D_B$ (so that the $A$- and $B$-particles perform the same random walk), then there exist constants $0 < C_i < \infty$ such that almost surely $\mathcal{C}(C_2 t) \subset B(t) \subset \mathcal{C}(C_1 t)$ for all large $t$, where $\mathcal{C}(r) = [-r, r]^d$. In a further paper, Kesten and Sidoravicius (2006), we proved a full "shape theorem" which states that $t^{-1}B(t)$ converges almost surely to a nonrandom compact set $B_0$ with the origin as an interior point, so the true growth rate for $B(t)$ is linear in $t$.

If $D_A \neq D_B$, then we could only prove the upper bound that $B(t) \subset \mathcal{C}(C_1 t)$ eventually. However, there is one extreme case for which a shape theorem and linear growth of $B(t)$ has also been proven. This is the so-called *frog model* in which $D_A = 0$, that is, the healthy particles stand still until they are infected [see Alves, Machado and Popov (2002) and Ramirez and Sidoravicius (2004)].

To get a better feel for the problem, we wanted to investigate the other extreme case, namely when $D_B = 0$. Taken literally, this is not an interesting case. In this case, the infected particles stand still and act as traps for the healthy particles. All that happens with any given $A$-particle is that it walks around until it coincides with one of the $B$-particles, after which it also stands still. The infected set $\widetilde{B}(t)$ equals $\widetilde{B}(0)$ at all $t \geq 0$ and the speed at which the infection spreads is 0. To obtain something interesting, we have to allow the $B$-particles to move, at least at some times. The simplest situation is the one-dimensional one, that is, when $d = 1$. We chose to let a $B$-particle move one unit to the right when an $A$-particle jumps on top of it. According to our rules, all $A$-particles which were one unit to the right of the $B$-particle are turned into $B$-particles at the time of this jump. This leads to the model described in the abstract.

The model described in the abstract is of further interest because it is a one-dimensional version of the celebrated DLA model of Witten and Sander (1981). In this model on $\mathbb{Z}^d$, one again has a growing aggregate $A(t) \subset \mathbb{Z}^d$ and one starts with $A(1) = \{\mathbf{0}\}$ = the origin. Usually, $t$ is taken to run through the integers and $A(t)$ has cardinality $t$. $A(t+1)$ is obtained from $A(t)$ by adding one point of $\mathbb{Z}^d$. This added point is the first point of the boundary of $A(t)$ which is reached by a random walker which starts at infinity [see Kesten (1987) for a more precise description]. The main difference between the model in the abstract and the DLA model of Witten and Sander is that the latter adds one $A$-particle to the system at a time, while in the former, there are infinitely many $A$-particles from the start. However, there have been various investigations for related models in which new $A$-particles are added to the system before all previously released $A$-particles have reached the boundary of the aggregate and are removed from the system; see, for instance, Lawler, Bramson and Griffeath (1992). In the physics literature, almost the same model as we discuss here was already studied by simulations in Voss (1984). However, in that paper, the $A$-particles do not perform independent random walks, but the system of $A$-particles evolves as an exclusion process; moreover, Voss (1984) considers the two-dimensional case. Also, Chayes and Swindle (1996) investigated hydrodynamic limits for the one-dimensional case in which the $A$-particles follow exclusion dynamics. We remark that the particle density in an exclusion process is necessarily at most 1. As we shall see, in our model, the case when the particle density $\mu$ is less than 1 can be handled much better than the case with $\mu \geq 1$. We have few results in the latter case.

As a side remark, we point out that DLA is usually considered in dimension $d > 1$, in which there is a whole new level of difficulty because we do not know how to describe the "shape" of $A(t)$.

Let us now turn to the problem raised in the abstract, namely the rate at which $R(t)$ grows. We take $\tau_0 = 0$. As stated, we take $R(0) = R(\tau_0) = 0$ and $N(i, 0)$, $i \geq 1$, an i.i.d. sequence of mean-$\mu$ Poisson random variables. All particles perform independent continuous-time simple random walks with jump-rate $D$ until they are absorbed by the aggregate. Unless otherwise stated, by "simple random walk," we mean a symmetric simple random walk. It is convenient to let the particles continue as a simple random walk, even after absorption, by giving the particles also a color, white or black. We start with all particles white, but absorption of the particle by the aggregate is now represented by changing the color of the particle from white to black at the time of its absorption. However, the particle's path is not influenced by its color. After a particle turns black, it continues with a continuous-time simple random walk path. A black particle has no interaction with any other particle, nor does it influence the motion of $R(\cdot)$. Thus, $R$ is not increased at a time $t$ when a black particle jumps to $R(t)$. In the sequel, we shall always use this description of the system with colored particles.

$N(i, t)$ denotes the number of white particles at the space–time point $(i, t)$. We successively define stopping times $\tau_k$ and take $R(t) = k$ on the time interval $[\tau_k, \tau_{k+1})$. Moreover, it will follow by induction on $k$ that

$$\text{at time } \tau_k, \text{ there are no white particles in } [0, R(\tau_k)] = [0, k]. \tag{1.1}$$

We take $\tau_0 = \mathbf{0}$. If $\tau_k$ and the $N(i, \tau_k)$ have been determined, and $R(\tau_k) = k$ and (1.1) holds, then we take

$$\tau_{k+1} = \inf\{t > \tau_k : \text{some white particle jumps to position } R(\tau_k) = k\}. \tag{1.2}$$

Since the particles perform simple random walk and (1.1) holds, only white particles at position $k+1$ at time $\tau_{k+1}-$ can jump to $k$ at time $\tau_{k+1}$. If such a jump occurs, we take $R(\tau_{k+1}) = k+1$ [i.e., $R(\cdot)$ jumps up by 1 at time $\tau_{k+1}$] and change to black the color of all white particles which were at $R(\tau_k)+1 = k+1$ at time $\tau_{k+1}-$ (this includes the particle which jumped to $k$ at $\tau_{k+1}$). It is clear that then (1.1) with $k$ replaced by $k+1$ holds so that we can now define $\tau_{k+2}$, etc. It also follows from this description that

$$R(t) = k \qquad \text{for } \tau_k \le t < \tau_{k+1}. \tag{1.3}$$

REMARK 1. We briefly indicate in this remark how our process can be constructed as a Markov process with the strong Markov property. However, anyone willing to accept the strong Markov property without proof will want to skip such a construction.

As our sample space, we take

$$\Omega := \prod_{i=1}^{\infty} D([0,\infty), \mathbb{Z}_+),$$

that is, the countable product of cadlag paths from $[0,\infty)$ to the nonnegative integers. All our random variables are functions on $\Omega$. We start with a countably infinite number of particles, which we order in some way as $\rho_1, \rho_2, \ldots.$ At the sample point $(\omega_1, \omega_2, \ldots) \in \Omega$, the $i$th coordinate, $\omega_i$, is the path of the particle $\rho_i$. The starting positions, $\omega_i(0)$, of the various particles are specified by the initial point of our process and the displacements $\{\omega_i(t) - \omega_i(0)\}_{t \ge 0}$, $i = 1, 2, \ldots,$ are i.i.d. simple continuous-time random walk paths. Thus, a sample point specifies the positions of all particles at all times. The colors of all particles at any given time $t$ and $R(t)$ are then also determined, but we do not attempt to write down an explicit expression for these random variables. If $Y(t)$ is the state of our process at time $t$, then $Y(t)$ is a point of $\Sigma$, which is the collection of all sequences $\{r, (n_i, \eta_i), i \ge 1\}$ with $r, n_i \in \mathbb{Z}_+$ and $\eta_i \in \{W, B\}$. $Y(t) = \{r, (n_i, \eta_i), i \ge 1\}$ means that $R(t) = r$ and the position and color of $\rho_i$ are $n_i$ and $\eta_i$, respectively. Of course, the process of the i.i.d. paths of the particles $\rho_i, 1 \le i < \infty$, is a Markov process and this makes $\{Y_t\}$ also into a process with the simple Markov property. However, we have to allow the possibility of explosion; we must add a cemetery point $\partial$ to our state space to define $Y(t)$ as a Markov process for all time $t$. We do not know whether this alone will make $\{Y(t)\}$ into a strong Markov process which can start at each point in $\Sigma$. We shall therefore choose a smaller state space than $\Sigma$.

Explosion can happen in two ways. First, $\tau_\infty := \lim_{k\to\infty} \tau_k$ may be finite. Second, it may be that $\tau_{k+1} = \tau_k$. This happens if and only if at some time $t$, there are infinitely many white particles at $R(t)+1$. We do not want to continue our process after such a time. In fact, we shall not continue our process beyond the time

$$\widehat{\tau} = \inf\{t\colon \text{there are infinitely many particles of any color at some site } z\}.$$

Let $P^\sigma$ denote the measure governing the process $\{Y(t)\}$ conditioned to start at $\sigma$. We then choose as our state space for $\{Y(t)\}$ the set

$$\Sigma_0 := \{\sigma \in \Sigma\colon P^\sigma\{\widehat{\tau} \wedge \tau_\infty < \infty\} = 0\}.$$

This description of the state space is rather indirect, but one can now prove that if the process $\{Y(t)\}$ starts at a point $\sigma \in \Sigma_0$, then it does not explode and stays in $\Sigma_0$ for all time a.s. $[P^\sigma]$. Moreover, the restricted process has the strong Markov property. Finally, $\Sigma_0$ is nonempty. If the starting point $\sigma$ is chosen by taking

$$\begin{aligned} &R(0) = \mathbf{0} \text{ and all particles initially white and } N(i,0), i \geq 1,\\ &\qquad \text{as i.i.d. Poisson random variables with mean } \mu, \end{aligned} \tag{1.4}$$

then $\sigma$ lies a.s. in $\Sigma_0$.

We shall not prove any of these statements here. Proofs can be given in the same manner (but actually simpler) as in Section 2 of Kesten and Sidoravicius (2003b). The principal step which makes this proof work is showing that for any $\sigma \in \Sigma_0$, any $L, T \geq 0$ and $\varepsilon > 0$, one can find a $K \geq L$ such that

$$P^\sigma\{\text{some particle which starts in } [K+1,\infty) \text{ enters } [0,L] \text{ during } [0,T]\} \leq \varepsilon.$$

The same is true if the initial state is chosen as in (1.4).

Let us now state our results. Throughout, $\mathbf{0}$ denotes the origin and $\{S(t)\}_{t\geq 0}$ is a continuous-time simple symmetric random walk on $\mathbb{Z}$ with jump rate $D$. Unless otherwise stated, $S(0) = \mathbf{0}$. We use $P\{A\}$ for the probability of the event $A$ in various probability models and $E$ for expectation with respect to $P$. It should be unambiguous from the context which probability measure we are discussing. $C_i$ will denote a constant with value in $(0,\infty)$. Its value may vary from formula to formula. Our first theorem states that for any value of $\mu$, the common expectation of the $N(k,0)$, it is the case that

$$\limsup_{t\to\infty} \frac{1}{t} R(t) < \infty \qquad \text{a.s.} \tag{1.5}$$

THEOREM 1. *Assume that $R(0) = \mathbf{0}$ and that the $N(i,0), i \geq 1$, are i.i.d. mean-$\mu$ Poisson variables. Then (1.5) holds. In fact, there exist constants $0 < C_i < \infty$ such that*

$$P\{R(t) > C_1 t\} \leq C_2 \exp[-C_3 t]. \tag{1.6}$$

REMARK 2. Theorem 1 remains valid if the particles perform an asymmetric simple random walk, that is, each jump of the random walk is $+1$ or $-1$ with probability $p_+$ and $p_- = 1 - p_+$, respectively. No change in the proof is required for this more general case.

In view of Theorem 1, it is reasonable to conjecture that $\lim_{t\to\infty}(1/t)R(t)$ exists and is constant a.s. One might even assume that this limit is strictly positive, but a quick (and quite general) argument in the next theorem shows that if $\mu < 1$, then "there are not enough particles around" to make $R(t)$ grow linearly with time.

THEOREM 2. *Assume that* $\{N(i,0)\}_{i\ge 1}$ *is a stationary ergodic sequence and that* $E\{N(i,0)\} = \mu$. *If* $0 < \mu < 1$, *then*

$$\lim_{t\to\infty} \frac{R(t)}{(\log t)^2\sqrt{t}} = 0 \qquad a.s. \tag{1.7}$$

*Moreover,* $R(t)/\sqrt{t}$, $t \ge 1$, *is a tight family, that is,*

$$P\{R(t) \ge x\sqrt{t}\} \to 0 \qquad \text{as } x \to \infty, \text{ uniformly in } t \ge 1. \tag{1.8}$$

*If we assume that the initial conditions satisfy (1.4), then (1.7) can be strengthened to*

$$\limsup_{t\to\infty} \frac{R(t)}{\sqrt{t}} < \infty \qquad a.s. \tag{1.9}$$

REMARK 3. One can formulate a $d$-dimensional analog of our model and of Theorem 2. In this version, one works on $\mathbb{Z}^d$ and at time 0, the aggregate consists of the origin only, while at the site $x \ne \mathbf{0}$, there are $N(x,0)$ particles, with the $N(x,0)$, $x \in \mathbb{Z}^d \setminus \{\mathbf{0}\}$ i.i.d. Poisson variables of mean $\mu$. Again, all particles perform independent continuous-time simple random walks. They all start out as white particles. We denote the aggregate at time $t$ by $A(t)$. If, at some time $t$, a white particle jumps from a site $x \notin A(t-)$ onto the aggregate, then we set $A(t) = A(t-) \cup \{x\}$ and all particles which were at $x$ at time $t-$ are changed to black at time $t$.

We define an outer radius of the aggregate by

$$R^{(o,d)} := \sup\{\|x\|_2 : x \in A(t)\}$$

and an inner radius as

$$R^{(i,d)}(t) := \inf\{\|x\|_2 : x \notin A(t)\}.$$

The latter is the distance from the origin to the nearest vertex outside $A(t)$. For this model Theorem 1 remains valid. More precisely, (1.6) and (1.5) with

$R(t)$ replaced by $R^{(o,d)}(t)$ still hold. Theorem 2 has the following analogue: if $\mu < 1$, then

$$\limsup_{t\to\infty} \frac{R^{(i,d)}(t)}{\sqrt{t}} < \infty \qquad \text{a.s.} \tag{1.10}$$

[Note that (1.10) is trivially true if there exists a site $x_0$ which is never occupied by $A(t)$.] We shall not give the proofs of these results here. They are essentially the same as for Theorem 1 and for (1.9).

If we strengthen our assumptions on the $N(i,0)$, then we can show that in the one-dimensional model, $R(t)/\sqrt{t}$ is actually bounded away from 0 in distribution. This holds for all $\mu > 0$.

THEOREM 3. *Assume that the $N(i,0)$, $i \geq 1$, are i.i.d. with finite second moment $\mu_2 > 0$. Then, for all $\varepsilon > 0$, there exists an $\eta = \eta(\varepsilon) > 0$ and a $t_0 = t_0(\varepsilon)$ such that*

$$P\left\{\frac{R(t)}{\sqrt{t}} > \eta\right\} \geq 1 - \varepsilon \qquad \textit{for all } t \geq t_0. \tag{1.11}$$

Unfortunately, the simple proof of (1.7) breaks down when $\mu > 1$ and we therefore conjecture that there exists a critical value $\mu_c \geq 1$ such that

$$\lim_{t\to\infty} \frac{1}{t} R(t) \text{ exists and is a.s. a constant which is } \begin{cases} > 0, & \text{if } \mu > \mu_c, \\ = 0, & \text{if } \mu < \mu_c. \end{cases} \tag{1.12}$$

A stronger conjecture would be that

$$\mu_c = 1. \tag{1.13}$$

Simulations certainly indicate that this is the case; see Figures 1–3 which plot $\log[(1/n)\sum_{i=1}^n R_i(t)]$ for various values of $n, t$, as a function of $\log t$, where $R_1, \ldots, R_n$ are independent copies of $R(t)$. We nevertheless marked the vertical axis as $\log ER(t)$ because we regard $(1/n)\sum_{i=1}^n R_i(t)$ as an approximation of $ER(t)$.

We have made only little progress toward proving (1.12), so we pose this as a problem.

OPEN PROBLEM 1. Prove (1.12) and, if this holds, determine $\mu_c$. If one becomes even more ambitious, one can ask whether power laws exist as $\mu \downarrow \mu_c$ and what the critical exponents are. To formulate this problem, we have to assume that $\lim_{t\to\infty}(1/t)R(t)$ exists. Let us write $Z(\mu)$ for this limit.

OPEN PROBLEM 2. Does

$$\lim_{\mu\downarrow\mu_c} \frac{\log Z(\mu)}{\log(\mu - \mu_c)}$$

exist and if so, what is its value?

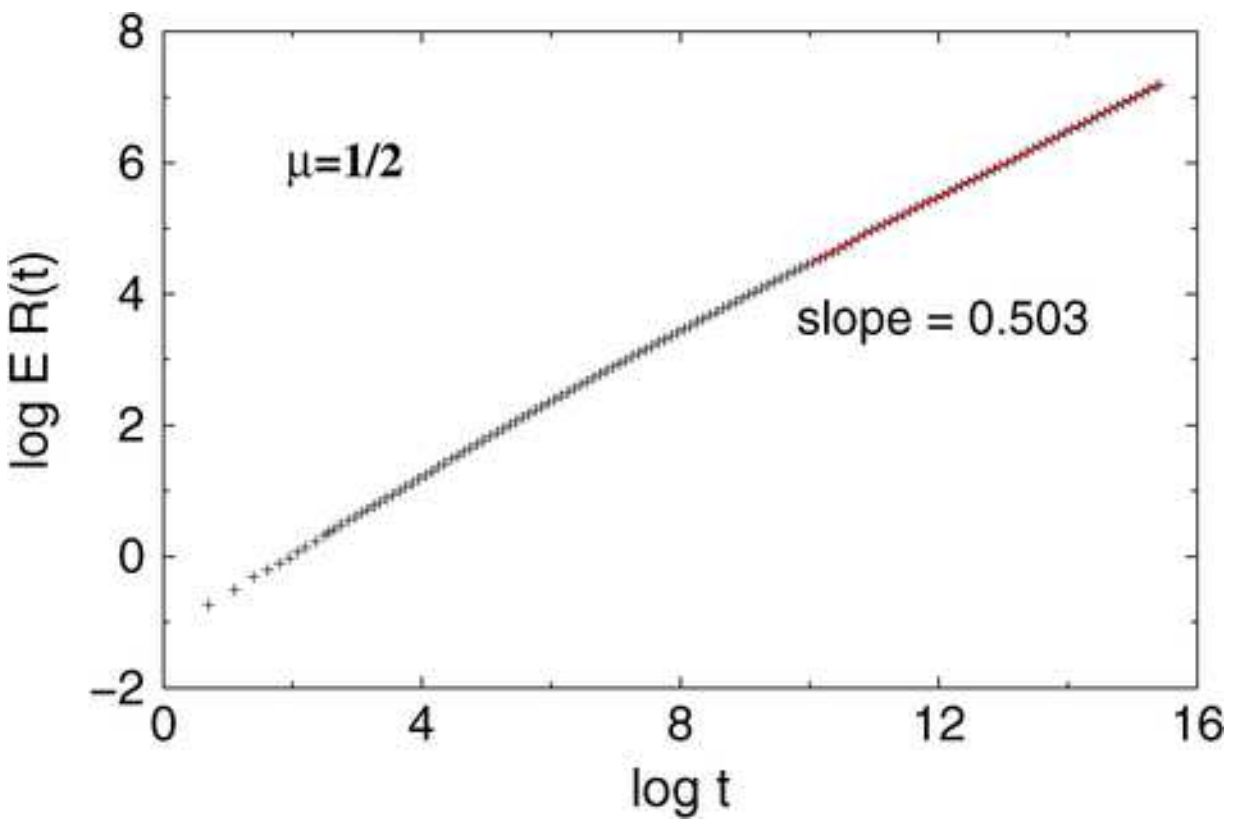

Fig. 1. *Graph of* $\log[(1/n)\sum_{i=1}^n R_i(t)]$ *against* $\log t$ *when* $\mu = 0.5$ *and* $n = 1000$ *(at least for part of the graph). The* $R_i(\cdot)$ *are independent runs of the process. The slope of the regression line is 0.503. The theory predicts a slope of* 0.5.

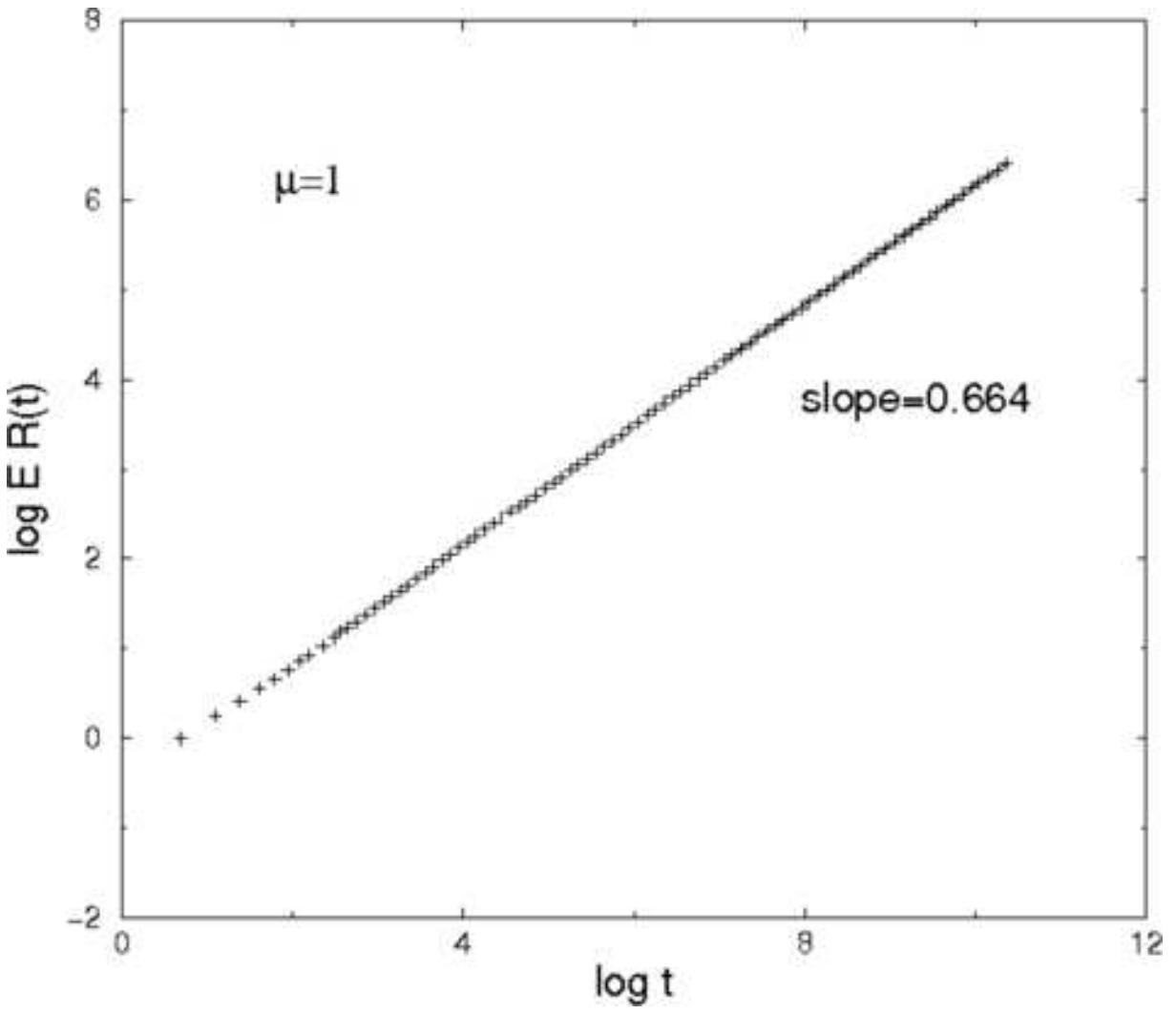

Fig. 2. *Graph of* $\log[(1/n)\sum_{i=1}^n R_i(t)]$ *against* $\log t$ *when* $\mu = 1$ *and* $n = 1000$. *The slope of the regression line is 0.664.*

A final problem about the DLA model is motivated by Theorems 2 and 3.

Open problem 3. Does $t^{-1/2}R(t)$ have a limit distribution as $t \to \infty$ when $\mu < 1$? It has been suggested to us that this problem could perhaps be handled by means of the techniques for establishing a hydrodynamic limit result for $R(t)$. Because of ignorance, we have made only a weak and unsuccessful effort in this direction.

The obvious approach to proving that $R(t)$ grows linearly in $t$ is to study our system as seen from the right edge of the aggregate. Indeed, the collection of positions of the white particles relative to $R(t)$ forms a Markov process. Does this Markov process have a nontrivial invariant probability distribution and, if so, is the invariant distribution unique? (By "nontrivial," we mean that we exclude the distribution which puts no particles at all to the right of the aggregate.) On an intuitive level, one would like to say that the invariant measure puts at position $R(t)+x$ roughly a Poisson number of particles with mean equal to $\mu$ times the probability that a particle at $R(t)+x$ is white. That is, the mean number of particles at $R(t)+x$ should be $\lim_{t\to\infty} \mu\nu(x,t)$, where

$$\nu(x,t) = P\{R(t) + x - S(s) > R(t-s) \text{ for } 0 \le s \le t\}.$$

Actually, all we want to know in first instance is that the density of white particles directly in front of $R(t)$ is bounded away from 0 as $t \to \infty$. We want to show that the system does not develop large holes without white particles in front of $R(t)$. To obtain such a result, we need some a priori control of $R(t) - R(t-s)$, which we do not know how to control. T. Kurtz (private communication) showed us that, conditionally on the $\sigma$-field generated by $\{R(s) : s \le t\}$, the $N(R(t)+x,t)$ have a Poisson distribution with a mean $\mu(x,t)$, say, even derived a system of differential equations for the $\mu(x,\cdot)$. Unfortunately, this system still involves the unknown random function $R(\cdot)$ in boundary conditions and we have been unable to make use of these differential equations.

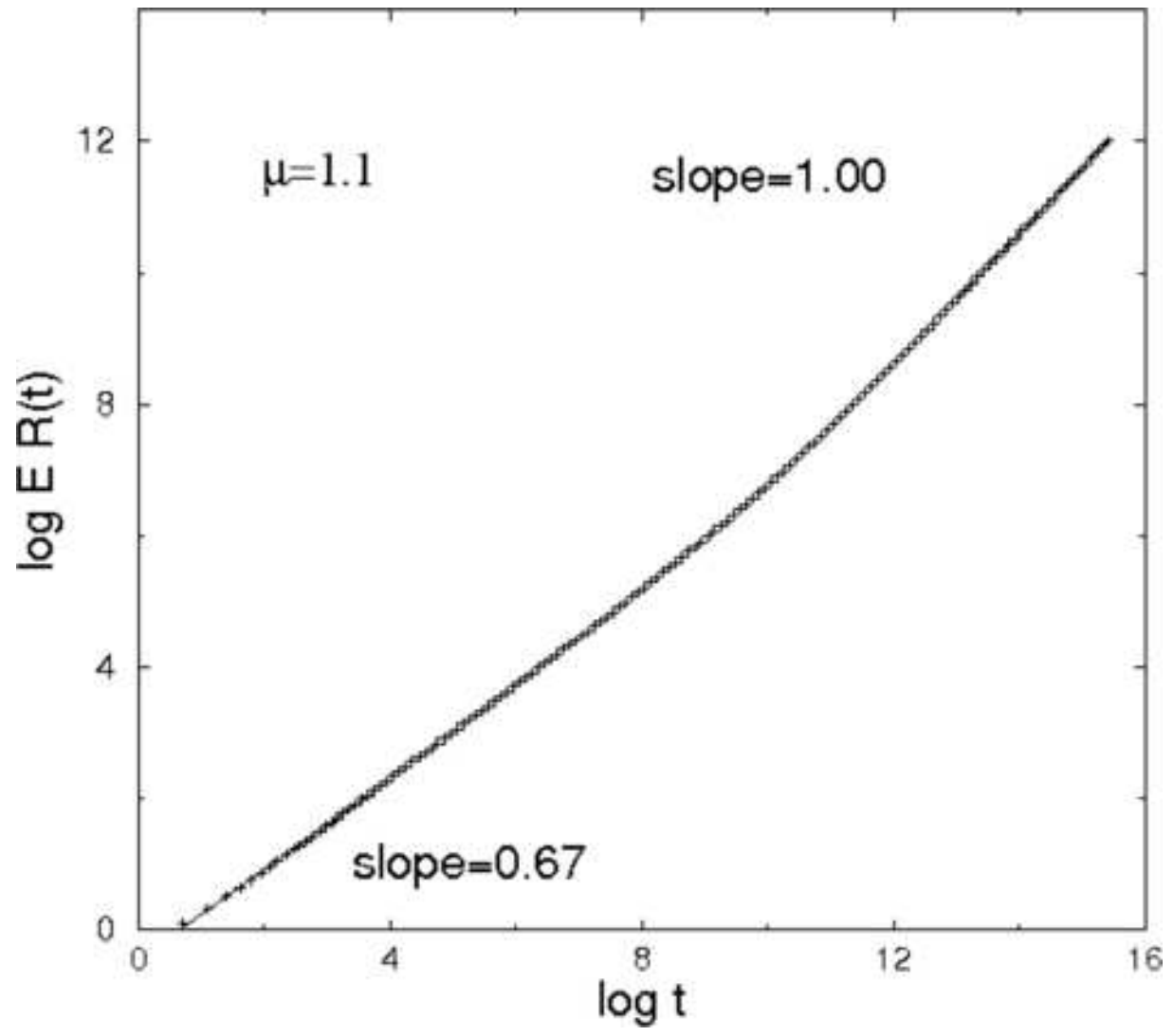

FIG. 3. *Graph of* $\log[(1/n)\sum_{i=1}^{n} R_i(t)]$ *against* $\log t$ *when* $\mu = 1.1$ *and* $n = 100$. *The slope of the fitted curve approaches 0.999 as t becomes large.*

Since we were unsuccessful in proving the existence of a nontrivial invariant probability measure for the Markov process of the last paragraph, we designed some caricatures of the model. We hope that these caricatures can be regarded as "approximations" to the true model and will help us to treat the true model. These caricatures have built-in mechanisms that make it more difficult for a large hole to form in front of the aggregate.

Caricature I. In this version, we still have an aggregate $A(t) = [0, R(t)] \cap \mathbb{Z}$. At time 0, we again put $N(i,0)$ white particles at $i$, with $\{N(i,0)\}$, $i \geq 1$, an i.i.d. sequence of mean-$\mu$ Poisson variables. In addition, we add $J$ particles to the system at some deterministic positions $x_1 \geq 1, \dots, x_J \geq 1$. Again, the right edge of the aggregate [i.e., $R(t)$] will increase at the successive times $\tau_k$ at which a particle jumps from $R(\tau_k-)+1$ to $R(\tau_k-)$. As before, at such a time, $R(\cdot)$ increases by 1 and changes to black the color of all the particles at $R(\tau_k-)+1$. Equation (1.3) again holds in this caricature. The difference between this model and the true one is in the motion of the particles, or rather in the time at which the particles start moving. This can be described by introducing another color. At the start, only the $J$ additional particles placed at $x_1, \dots, x_J$ will be white. The other particles [$N(i,0)$ of them at position $i$] will be colored red. Red particles do not move. Once a particle turns white or black, it performs a continuous-time simple random walk, as in the true model. These random walks are independent of each other. A particle changes from red or white to black when it is first at position $k$ at time $\tau_k-$ for some $k$. At all times, there will be exactly $J$ white particles in the system. If, at time $\tau_k$, $m$ white particles change to black, then we replenish the system by changing $m$ red particles to white, which then begin their random walks. To complete the description, we have to specify how the $m$ red particles which become white are chosen. We will pick these as close to the right edge as possible. That means that at time $\tau_k$, we first change particles at $k+1$ from red to white. If there are at least $m$ red particles at $k+1$ at time $\tau_k-$, then we change exactly $m$ of these to white and no other red particles turn white at this moment. If there are $m_1 < m$ red particles at $k+1$ at time $\tau_k-$, then we change all $m_1$ of these red particles to white. We then look for $m-m_1$ red particles at $k+2$. As before, if there are $m_2$ red particles at $k+2$ and $m_2 \geq m-m_1$, then we turn exactly $m-m_1$ red particles at $k+2$ into white and do not change any further red particles to white. If $m_2 < m-m_1$, then we change all $m_2$ of the red particles at $k+2$ to white and we still need to change $m-m_1-m_2$ particles. We now search for these at $k+3$, etc., until we have changed $m$ red particles to white.

This version indeed has some of the desirable properties. Equation (1.5) still holds for this model. In fact, since, at any time, there are only $J$ white particles present in the system, the rate at which $R(\cdot)$ jumps is at most $JD/2$, so

$$\limsup_{t\to\infty} \frac{1}{t} R(t) \le \frac{JD}{2} \qquad \text{a.s.} \tag{1.14}$$

Furthermore, if $0 < \mu < 1$, then (1.7) and (1.8) hold for this caricature. No essential changes in the proof of Theorem 2 are needed for these. It is, however, not clear whether (1.11) holds when $\mu < 1$. Most importantly, we can show in this model that there exists a $J_0$ such that if $J \ge J_0$ and $\mu \ge 2J$, then

$$\lim_{t\to\infty} \frac{1}{t} R(t) \text{ exists and is strictly positive a.s.} \tag{1.15}$$

Thus, there is a phase transition in the large-time growth rate of $R(t)$ in this model.

Unfortunately, the proof of (1.15) for this caricature is still rather complicated. Since this is only a caricature, we shall not give this proof, but instead treat a simpler caricature, one which is a bit further removed from the true model.

Caricature II. This caricature corresponds more or less to Caricature I with $\mu = \infty$. That is, we start with infinitely many red particles at each $i \ge 1$, plus $J$ additional white particles at $x_1, \dots, x_J$. Everything but the choice of which red particles to turn into white ones is as in Caricature I. However, if at time $\tau_k$, $m$ red particles have to be turned into white ones, then we select these at positions $k+Y_1, k+Y_2, \dots, k+Y_m$ with the $Y_1, \dots, Y_m$ i.i.d. with some common distribution $G$, where $G$ is concentrated on $\{1, 2, \dots\}$. In Section 4, we shall show that if $J$ is large and $G$ has a suitable number of moments, then (1.15) again holds for this caricature (see Theorem 4). Also, (1.1) still holds in this model.

The reason why this caricature is relatively simple to treat is that the positions of the $J$ white particles (as seen from the front of the aggregate) form a Markov process with a countable state space. There is a standard method to prove positive recurrence of such processes, namely to apply Foster's criterion after finding a suitable positive supermartingale or Lyapounov function [see Fayolle, Malyshev and Menshikov (1995), Section 2.2]. This is indeed the method which we shall use in Section 4.

## 2. A linear upper bound for $R(t)$.

PROOF OF THEOREM 1. This proof is actually contained in the proof of Proposition 4 and Theorem 1 of Kesten and Sidoravicius (2005). However, in the special case with which we deal here, the proof simplifies and we can quickly show the principal step. Basically, this is a Peierls argument, in that it estimates the expected number of certain paths. For the proof of Theorem 1, it is convenient to label the particles in a different way than in Section 1. We shall order the particles at $i$ at time 0 in an arbitrary way and denote the $j$th particle at $i$ at time 0 by $\langle i, j\rangle$. We say that $\langle i, j\rangle$ *exists* if $N(i, 0) \geq j$. Then, if we take $k = \lceil C_1 t\rceil$,

$$\begin{aligned} P\{R(t) \geq C_1 t\} &\leq P\{\tau_k \leq t\} \\ &\leq \int\cdots\int_{0\leq t_1\leq t_2\leq\cdots\leq t_k\leq t} P\{\text{some existing particle } \langle u_i, v_i\rangle \text{ jumps} \\ &\qquad\qquad \text{from } i \text{ to } i-1 \text{ during } dt_i, 1 \leq i \leq k, \\ &\qquad\qquad \text{and } \langle u_i, v_i\rangle, 1 \leq i \leq k, \text{ are distinct}\}. \end{aligned} \tag{2.1}$$

Note that the $\langle u_i, v_i\rangle$ have to be distinct because a given particle can change from white to black at most once. As in (2.51) of Kesten and Sidoravicius (2005), the right-hand side here equals

$$\begin{aligned} &\int\cdots\int_{0\leq t_1\leq t_2\leq\cdots\leq t_k\leq t} [D/2]^k\, dt_1\cdots dt_k \\ &\qquad \times \sum_{\langle u_i,v_i\rangle \text{ distinct}} E\prod_{i=1}^{k}[I[\langle u_i, v_i\rangle \text{ exists}]P\{S(t_i-) = i+1-u_i\}] \\ &\leq \int\cdots\int_{0\leq t_1\leq t_2\leq\cdots\leq t_k\leq t} dt_1\cdots dt_k\left[\frac{D\mu}{2}\right]^k = \left[\frac{D\mu t}{2}\right]^k\frac{1}{k!} \\ &\leq C_2\left[\frac{De\mu t}{2k}\right]^k. \end{aligned} \tag{2.2}$$

The bound (1.6) with $C_1 = De\mu$ now follows (recall that $k = \lceil C_1 t\rceil$) and $C_3 = C_1 \log 2$. □

We already stated in Remark 1 that $\min(\widehat{\tau}, \tau_\infty)$ is almost surely infinite. Theorem 1 confirms that [under (1.4)] $\tau_\infty = \infty$. The next lemma confirms that also $\widehat{\tau} = \infty$. We show the short proof [which is also part of the proof of Proposition 4 in Kesten and Sidoravicius (2005)] because some of the computations in the proof will be needed again in the next section. Define

$$\alpha_s(z) = P\{S(s) = -z\}.$$

Note that

$$\alpha_{s+u}(z) \geq e^{-Du}\alpha_s(z) \quad \text{and} \quad \alpha_{s+u}(z) \geq \alpha_u(z')\alpha_s(z-z') \tag{2.3}$$

for any $z, z' \in \mathbb{Z}$ and that

$$\begin{aligned} &\int_0^{t+1} \alpha_s(y-z)\,ds \\ &\quad = E\{\text{amount of time spent by } S(\cdot) \text{ at } z \text{ during } [0,t+1] | S(0)=y\} \\ &\quad \geq P\{y+S(s) \text{ reaches } z \text{ at some } s \leq t \\ &\qquad\quad \text{and stays at } z \text{ at least one unit of time}\} \qquad (2.4) \\ &\quad \geq e^{-D} P\{y+S(s)=z \text{ for some } s \leq t\}. \end{aligned}$$

Lemma 1. *Assume that the* $N(i,0)$, $i \geq 1$, *are i.i.d. mean-*$\mu$ *Poisson variables. Then* $\widehat{\tau} = \infty$ *almost surely.*

Proof. For any $t \geq 0, z \in \mathbb{Z}$,

$$\begin{aligned} E\Big\{\sup_{s \leq t} N(z,s)\Big\} &\leq E\{\text{number of particles which visit } z \text{ during } [0,t]\} \\ &\leq \sum_{y\in\mathbb{Z}} \mu P\{y+S(s)=z \text{ for some } s \leq t\} \\ &\leq \sum_{y\in\mathbb{Z}} \mu e^D \int_0^{t+1} \alpha_s(y-z)\,ds \qquad [\text{by } (2.4)] \\ &= \mu e^D(t+1) < \infty. \end{aligned}$$

Thus, $P\{N(z,s)=\infty \text{ for some } z \in \mathbb{Z} \text{ and } s < \infty\} = 0$. □

**3. A sublinear upper bound for $R(t)$ when $\mu < 1$.** In this section, we show that $R(t)$ cannot grow linearly with $t$ when $\mu < 1$. This results requires far less than (1.4), as shown in Theorem 2. If (1.4) is assumed and $\mu < 1$, then Theorems 2 and 3 show that $R(t)$ is of order $\sqrt{t}$.

Proof of Theorem 2. Define

$$U(t) = (\text{number of black particles in the system at time } t)$$

and

$$V(r,t) = (\text{number of particles which moved into } [0,r] \text{ during } [0,t]).$$

Note that $V(r,t)$ only counts particles which were outside $[0,r]$ at time 0 and that a particle does not have to be in $[0,r]$ *at* time $t$ to be counted in $V(r,t)$. A particle can be black at time $t$ only if it started in $[0,r]$ or if it coincided with $R(s)$ at some time $s \le t$. Therefore, one has, on the event $\{R(t) = r\}$,

$$\begin{aligned} U(t) &\le \sum_{i=0}^{r} N(i,0) \\ &\quad + (\text{number of particles which moved into } [0,r] \text{ during } [0,t]) \\ &= \sum_{i=0}^{r} N(i,0) + V(r,t). \end{aligned} \tag{3.1}$$

On the other hand, we must have

$$U(t) \ge r \tag{3.2}$$

since at least one new particle turns black at each time when $R(t)$ increases by 1. Thus, still on $\{R(t) = r\}$,

$$r \le \sum_{i=1}^{r} N(i,0) + V(r,t). \tag{3.3}$$

Now, if the $N(i,0)$ form a stationary ergodic sequence with common mean $\mu$, then

$$\sum_{i=0}^{r} N(i,0) = \mu r + o(r) \le [\mu + (1-\mu)/2]r$$

for $r \ge$ some (random) $r_0$ almost surely. Thus, for such $r$,

$$\tfrac{1}{2}(1-\mu)r \le V(r,t).$$

Now, the event $\{R(t) \ge A\}$ (with $A$ a positive integer) can occur only if, for some $s \le t$, $R(s) = A$. Therefore, for such an $s$,

$$\begin{aligned} \{R(t) \ge A\} &\subset \{A \le r_0\} \cup \left\{A \le \frac{2}{(1-\mu)} V(A,s)\right\} \\ &\subset \{A \le r_0\} \cup \left\{A \le \frac{2}{(1-\mu)} V(A,t)\right\}. \end{aligned} \tag{3.4}$$

Now, let us estimate $E\{V(A,t)\}$. This is the expected number of particles which start in $[A+1,\infty)$ and which enter $[0,A]$ during $[0,t]$. Thus, for suitable constants $0 < C_i < \infty$, independent of $A$, and $t \ge 1$,

$$E\{V(A,t)\} = \sum_{i \ge A+1} E\{N(i,0)\} P\left\{\inf_{s \le t}[i + S(s)] \le A\right\}$$

$$
\begin{aligned}
&= \mu \sum_{i=A+1}^{\infty} P\Big\{\inf_{s\le t} S(s) \le A-i\Big\} \\
&\le \mu \sum_{\ell=1}^{\infty} C_1 \exp\bigg[-C_2 \frac{\ell^2}{t+\ell}\bigg] \qquad (3.5)\\
&\qquad \text{[e.g., by the inequality (2.42) in Kesten and Sidoravicius (2003a)]}\\
&\le C_3\sqrt{t}.
\end{aligned}
$$

Consequently,

$$
P\bigg\{A \le \frac{2}{(1-\mu)} V(A,t)\bigg\} \le C_4 \frac{\sqrt{t}}{A}. \qquad (3.6)
$$

In particular, we obtain for $t = 4^k$ and $A = \lceil \varepsilon k^2 2^k \rceil$, that for any fixed $\varepsilon > 0$,

$$
P\bigg\{\lceil \varepsilon k^2 2^k \rceil \le \frac{2}{(1-\mu)} V(\lceil \varepsilon k^2 2^k \rceil, 4^k)\bigg\} \le \frac{C_4}{\varepsilon k^2}.
$$

Thus, by Borel–Cantelli, almost surely $\frac{2}{(1-\mu)} V(\varepsilon k^2 2^k, 4^k) \le \lceil \varepsilon k^2 2^k \rceil$ for all large $k$. Also, $\varepsilon k^2 2^k > r_0$ for all large $k$ almost surely. (3.4) now tells us that almost surely,

$$
R(4^k) \le \lceil \varepsilon k^2 2^k \rceil, \qquad \text{eventually}.
$$

Since $R(\cdot)$ is nondecreasing, this implies (1.7).

The tightness in (1.8) follows in a similar way from (3.4) and (3.6).

If (1.4) holds, then $V(A,t)$ is actually bounded by a Poisson random variable with mean at most

$$
\sum_{i \ge A+1} E\{N(i,0)\} P\Big\{\inf_{s\le t}[i + S(s)] \le A\Big\} = \mu \sum_{i\ge 1} P\Big\{\inf_{s\le t} S(s) \le -i\Big\} \le C_3\sqrt{t}.
$$

If we then take $t = 4^k$ and $A = A_k =: [4C_3/(1-\mu)]2^k$, we get, instead of (3.6),

$$
P\bigg\{V(A_k, 4^k) \ge \frac{(1-\mu)A_k}{2}\bigg\} \le \exp\bigg[-\frac{\theta}{2}(1-\mu)A_k + (e^\theta - 1)C_3 2^k\bigg]
$$

for any $\theta \ge 0$. In particular, if $\theta > 0$ is taken such that $e^\theta - 1 \le 3\theta/2$, then we obtain

$$
P\bigg\{V(A_k, 4^k) \ge \frac{(1-\mu)A_k}{2}\bigg\} \le \exp\bigg[-\frac{\theta}{2} C_3 2^k\bigg].
$$

Thus, in this situation, a.s. $V(A_k, 4^k) < (1-\mu)A_k/2$ and $R(4^k) < A_k$, eventually, and (1.9) follows from the monotonicity of $R(\cdot)$ as before. $\square$

Proof of Theorem 3. We shall give a proof by contradiction. Let $\varepsilon > 0$ be given and assume that (1.11) fails for this $\varepsilon$. There then exists some sequence $t_n \to \infty$ such that for all $\eta > 0$,

$$P\left\{\frac{R(t_n)}{\sqrt{t_n}} > \eta\right\} < 1 - \varepsilon \qquad \text{for all large } n. \tag{3.7}$$

This, together with the monotonicity of $R(\cdot)$, implies that for all $0 < \gamma \le 1, \eta > 0$ and large enough $n$,

$$P\left\{\frac{R(\gamma t_n)}{\sqrt{t_n}} \le \eta\right\} \ge P\left\{\frac{R(t_n)}{\sqrt{t_n}} \le \eta\right\} \ge \varepsilon. \tag{3.8}$$

In order to choose $\gamma$ and $\eta$, we need some preparations. Fix $\alpha > 0$ such that

$$P\left\{\inf_{s \le t} S(s) \le -2\sqrt{t}\right\} \ge \alpha \qquad \text{for all } t \ge 1.$$

Such an $\alpha > 0$ exists by the central limit theorem. Define

$$W(t) = (\text{number of particles which start in } [\sqrt{t}, 2\sqrt{t}] \\ \text{at time 0 and reach the origin during } [0, t]).$$

Then each particle which starts in $[\sqrt{t}, 2\sqrt{t}]$ at time 0 has a probability of at least $\alpha$ of reaching the origin during $[0, t]$. Thus, conditionally on $\sum_{i \in [\sqrt{t}, 2\sqrt{t}]} N(i, 0) = N$, $W(t)$ is stochastically larger than $B(N, \alpha)$, where $B(N, \alpha)$ is a binomially distributed random variable corresponding to $N$ trials, each with success probability $\alpha$. In particular,

$$P\{W(t) \le \tfrac{1}{2}\alpha\mu\sqrt{t}\} \le P\left\{\sum_{i \in [\sqrt{t}, 2\sqrt{t}]} N(i, 0) \le \tfrac{3}{4}\mu\sqrt{t}\right\} \\ + P\{B(\tfrac{3}{4}\mu\sqrt{t}, \alpha) \le \tfrac{1}{2}\alpha\mu\sqrt{t}\}.$$

By simple weak law of large numbers estimates (i.e., Chebyshev's inequality), the right-hand side here tends to 0 as $t \to \infty$ so that for all large $n$,

$$P\{W(t_n) \le \tfrac{1}{2}\alpha\mu\sqrt{t_n}\} \le \tfrac{\varepsilon}{4}. \tag{3.9}$$

We also define

$$\widetilde{W}(\gamma, t) := (\text{number of particles which start in } [\sqrt{t}, 2\sqrt{t}] \\ \text{at time 0 and enter } (-\infty, \tfrac{1}{2}\sqrt{t}] \text{ during } [0, \gamma t])$$

and let $\beta$ be such that

$$P\left\{\inf_{s \le \gamma t} S(s) \le -\tfrac{1}{2}\sqrt{t}\right\} \le \beta \qquad \text{for all } t \ge 1.$$

Then

$$E\{\widetilde{W}(\gamma,t)\} \le \sum_{i\in[\sqrt{t},2\sqrt{t}]} EN(i,0)P\Big\{\inf_{s\le\gamma t} S(s) \le -\tfrac{1}{2}\sqrt{t}\Big\} \le (\mu\sqrt{t}+2)\beta,$$

where $\mu = EN(i,0)$. Note that we can take $\beta$ arbitrarily small by taking $\gamma$ small [by means of (2.42) in Kesten and Sidoravicius (2003a)]. We can therefore also fix $\gamma > 0$ and $\beta$ correspondingly small such that

$$P\{\widetilde{W}(\gamma,t_n) \ge \tfrac{1}{4}\alpha\mu\sqrt{t_n}\} \le \frac{\varepsilon}{4} \tag{3.10}$$

for all large $n$. With $\gamma$ and $\beta$ fixed in this way, we have from (3.8)–(3.10), for any fixed $\eta > 0$ and large $n$, that

$$\begin{aligned} P\bigg\{ &\frac{R(\gamma t_n)}{\sqrt{t_n}} \le \frac{R(t_n)}{\sqrt{t_n}} \le \eta, \\ &W(t_n) \ge \frac{1}{2}\alpha\mu\sqrt{t_n} \text{ and } \widetilde{W}(\gamma,t_n) \le \frac{1}{4}\alpha\mu\sqrt{t_n}\bigg\} \ge \frac{\varepsilon}{2}. \end{aligned} \tag{3.11}$$

For the remainder of this proof, $t$ will always be restricted to belong to the sequence $\{t_n\}$, even if we do not attach a subscript to $t$. We now define

$$\begin{aligned} \mathcal{C}_j = {} & \text{collection of particles which start in } [\sqrt{t},2\sqrt{t}] \text{ and which} \\ & \text{turn black at the } j\text{th jump of } R(\cdot) \text{ after time } \gamma t. \end{aligned}$$

We use $|B|$ to denote the cardinality of a collection $B$. We claim that on the event in the left-hand side of (3.11), we have

$$\sum_{j=1}^{R(t)-R(\gamma t)} |\mathcal{C}_j| \ge (1/4)\alpha\mu\sqrt{t}. \tag{3.12}$$

Note that the jumps of $R(\cdot)$ in the time interval $(\gamma t, t]$ are precisely the $j$th jump after $\gamma t$ for some $1 \le j \le R(t) - R(\gamma t)$. To see (3.12), note that any particle which reaches the origin during $[0,t]$ must have coincided with $R(s)$ for some $s \le t$ and must therefore be black at time $t$. In particular, this holds for the $W(t)$ particles which start in $[\sqrt{t},2\sqrt{t}]$ and which reach the origin during $[0,t]$. If we restrict ourselves to $\eta < 1/2$, then, on the event in the left-hand side of (3.11), it is the case that $R(s) \le R(\gamma t) \le \eta\sqrt{t} \le (1/2)\sqrt{t}$ for $s \le \gamma t$, so the particles which do not enter $(-\infty,(1/2)\sqrt{t}]$ during $[0,\gamma t]$ cannot have turned black yet at time $\gamma t$. Thus, on the event (3.11), we have at least $W(t) - \widetilde{W}(\gamma,t) \ge (1/4)\alpha\mu\sqrt{t}$ particles which start in $[\sqrt{t},2\sqrt{t}]$ and which turn black during $(\gamma t,t]$. All these particles belong to some $\mathcal{C}_\ell$ with $\gamma t < \tau_\ell \le t$ and are therefore counted in the left-hand side of (3.12) so that (3.12) follows.

Finally, we will show that we can choose $\delta > 0$ such that

$$P\left\{\sum_{j=1}^{R(t)-R(\gamma t)} |\mathcal{C}_j|^2 \le \delta\sqrt{t}\right\} \ge 1 - \frac{\varepsilon}{4}. \tag{3.13}$$

Then the probability that the events in the left-hand sides of (3.11) and (3.13) both occur is at least $\varepsilon/4$. However, on the intersection of these two events, we have, by (3.12) and Schwarz' inequality,

$$\tfrac{1}{16}\alpha^2\mu^2 t \le \left[\sum_{j=1}^{R(t)-R(\gamma t)} |\mathcal{C}_j|\right]^2 \le [R(t) - R(\gamma t)] \sum_{j=1}^{R(t)-R(\gamma t)} |\mathcal{C}_j|^2 \le R(t)\delta\sqrt{t}.$$

In particular, this implies that

$$R(t) \ge \frac{\alpha^2\mu^2}{16\delta}\sqrt{t}$$

on the intersection of the left-hand sides of (3.11) and (3.13). This, however, is impossible for $\eta < \alpha^2\mu^2/(16\delta)$ since one cannot simultaneously have $R(t) \le \eta\sqrt{t}$ and $R(t) \ge \alpha^2\mu^2\sqrt{t}/(16\delta)$. Thus, the assumption that (1.11) fails leads to a contradiction.

It remains to show (3.13). This will follow from a bound, on $E\{\sum_{j=1}^{R(t)-R(\gamma t)} |\mathcal{C}_j|^2\}$. Before we prove such a bound, we remind the reader of some basic inequalities. First, for some constant $C_1 < \infty$ depending on $D$ only, we have

$$\sup_k P\{S(s) = k\} \le \frac{C_1}{\sqrt{s}+1} \qquad \text{for } s > 0.$$

This follows from the local central limit theorem [see also the proof of Lemma 12 in Kesten and Sidoravicius (2003a) and Proposition 7.10 in Spitzer (1976)]. This estimate can be slightly refined. Indeed, for some further constants $0 < C_i < \infty$ depending on $D$ only, it holds uniformly in $k$ that

$$\begin{aligned}
P\{S(s) = k\} &= \sum_{\ell} P\{S(s/2) = \ell\}P\{S(s) - S(s/2) = k - \ell\} \\
&= \sum_{\ell \le |k|/2} P\{S(s/2) = \ell\}P\{S(s/2) = k - \ell\} \\
&\quad + \sum_{\ell > |k|/2} P\{S(s/2) = \ell\}P\{S(s/2) = k - \ell\} \\
&\le \frac{C_1}{\sqrt{s/2}+1} P\{S(s/2) \ge |k|/2\} + P\{S(s/2) > |k|/2\}\frac{C_1}{\sqrt{s/2}+1} \\
&\le \frac{C_2}{\sqrt{s}+1} \exp\left[-C_3 \frac{k^2}{s+|k|}\right],
\end{aligned}$$

where, in the last step, we used (2.42) of Kesten and Sidoravicius (2003a). In particular, if $\pi(\xi, s)$ denotes the position of a particle $\xi$ at time $s$, then for $k \geq 0$, $\gamma t \leq s \leq t$ and $z \in [\sqrt{t}, 2\sqrt{t}]$,

$$\begin{aligned} P\{\pi(\xi,s) = k+1 | \pi(\xi,0) = z\} &\leq \frac{C_4}{\sqrt{\gamma t}} \exp\left[-\frac{C_3(k+1-z)^2}{t+|k+1-z|}\right] \\ &\leq \frac{C_4}{\sqrt{\gamma t}} \exp\left[-\frac{C_5 k^2}{t+k}\right]. \end{aligned} \tag{3.14}$$

We turn to the proof of (3.13) proper. We start with the basic relation

$$\begin{aligned} \sum_{j=1}^{R(t)-R(\gamma t)} |\mathcal{C}_j|^2 &= \sum_{j=1}^{R(t)-R(\gamma t)} |\mathcal{C}_j| \sum_{\rho \in \mathcal{C}_j} 1 \\ &= \sum_{j=1}^{R(t)-R(\gamma t)} \sum_{\rho_1} \sum_{\rho_2} I[\rho_1 \text{ and } \rho_2 \text{ both turn black at} \\ &\qquad \text{time } \tau_{R(\gamma t)+j} \text{ when all} \\ &\qquad \text{particles of } \mathcal{C}_j \text{ turn black} \\ &\qquad \text{and this happens during } [\gamma t, t]] \\ &= \sum_{\rho_1} \sum_{\rho_2} I[\rho_1 \text{ and } \rho_2 \text{ turn black} \\ &\qquad \text{simultaneously at some time in } (\gamma t, t]]. \end{aligned} \tag{3.15}$$

Here, the sum over each $\rho_i$ is over all particles $\rho$ which start in $[\sqrt{t}, 2\sqrt{t}]$ and which change color during $(\gamma t, t]$. We denote the $\sigma$-field generated by the initial $\{N(i,0), i \geq 1\}$ by $\mathcal{F}_0$. We then have the following bound on the conditional expectation of (3.15), given the initial data:

$$\begin{aligned} &E\left\{\sum_{j=1}^{R(t)-R(\gamma t)} |\mathcal{C}_j|^2 | \mathcal{F}_0\right\} \\ &\quad \leq \sum_{\rho_1} \sum_{\rho_2} P\{\rho_1 \text{ and } \rho_2 \text{ change color at the same time in } (\gamma t, t] | \mathcal{F}_0\}. \end{aligned} \tag{3.16}$$

Here, the sum in the right-hand side is over all $\rho_1$ and $\rho_2$ which start in $[\sqrt{t}, 2\sqrt{t}]$. After taking the expectation over the initial state, the contribution to $E\{\sum_{j=1}^{R(t)-R(\gamma t)} |\mathcal{C}_j|^2\}$ from pairs with $\rho_1 = \rho_2$ is at most

$$E\left\{\sum_{\rho_1} 1\right\} = E\left\{\sum_{i \in [\sqrt{t}, 2\sqrt{t}]} N(i,0)\right\} \leq \mu(\sqrt{t}+2). \tag{3.17}$$

To handle pairs $\rho_1 \neq \rho_2$, define $R^r(s;\rho_1)$ and $R^r(s;\rho_1,\rho_2)$ to be the values of $R(s)$ in the system from which $\rho_1$, respectively $\rho_1$ and $\rho_2$, have been removed at time 0. Now $\rho_1$ and $\rho_2$ may change color simultaneously in $ds$ in three ways: (i) $\rho_2$ jumps during $ds$ from $R^r(s-;\rho_1)+1$ to $R^r(s-;\rho_1)$ and $\rho_1$ is at $R^r(s-;\rho_1)+1$ at time $s-$ and then changes color at the time when $\rho_2$ jumps and $R^r(\cdot;\rho_1)$ increases by 1; (ii) the scenario in (i) with the roles of $\rho_1$ and $\rho_2$ interchanged is followed; (iii) some white particle $\rho_3$, different from $\rho_1$ and $\rho_2$, jumps during $ds$ from $R^r(s-;\rho_1,\rho_2)+1$ to $R^r(s-;\rho_1,\rho_2)$, and $\rho_1$ and $\rho_2$ are at $R^r(s-;\rho_1,\rho_2)+1$ at time $s-$ and then change color at the time when $\rho_3$ jumps and $R^r(\cdot;\rho_1,\rho_2)$ increases by 1.

Now, observe that for $\rho_1 \neq \rho_2$, $R^r(\cdot;\rho_1)$ and $\pi(\rho_2,\cdot)$ are independent of the path of $\rho_1$. Therefore, the conditional probability of (i) taking place, given $\mathcal{F}_0$, and for a given $\rho_1,\rho_2$ with $\rho_2 \neq \rho_1$ which start in $[\sqrt{t}, 2\sqrt{t}]$, is at most

$$\begin{aligned}
&\sum_{k\geq 0}\int_{\gamma t}^{t} P\{R^r(s-;\rho_1)=k,\\
&\qquad\qquad \pi(\rho_2,s-)=k+1 \text{ and } \rho_2 \text{ is white at time } s-|\mathcal{F}_0\}\\
&\qquad\quad \times P\{\pi(\rho_1,s-)=k+1\}\frac{D}{2}\,ds\\
&\quad \leq \sum_{k\geq 0}\int_{\gamma t}^{t} P\{R^r(s-;\rho_1)=k,\\
&\qquad\qquad\qquad \pi(\rho_2,s-)=k+1 \text{ and } \rho_2 \text{ is white at time } s-|\mathcal{F}_0\}\\
&\qquad\qquad \times \frac{C_4 D}{2\sqrt{\gamma t}}\exp\left[-\frac{C_5 k^2}{t+k}\right]ds \qquad \text{[by (3.14)].}
\end{aligned} \tag{3.18}$$

We now sum this first over the $\rho_2 \neq \rho_1$ which start in $[\sqrt{t}, 2\sqrt{t}]$. For fixed $\rho_1$ and $k$, the events that $R$ jumps from $k$ to $k+1$ due to a jump of $\rho_2$ in $ds$ are disjoint for different $\rho_2$ and $s$. Therefore, the sum of (3.18) over $\rho_2$ and integral over $s$ is at most

$$\begin{aligned}
&\frac{C_4 D}{2\sqrt{\gamma t}}\sum_{k\geq 0}\int_{\gamma t}^{t} P\{R^r(\cdot;\rho_1) \text{ jumps from } k \text{ to } k+1 \text{ during } ds \text{ due to}\\
&\qquad\qquad\qquad\qquad \text{a jump of some particle other than } \rho_1\}\\
&\qquad \times \exp\left[-\frac{C_5 k^2}{t+k}\right]\\
&\quad \leq \frac{C_4 D}{2\sqrt{\gamma t}}\sum_{k\geq 0}\exp\left[-\frac{C_5 k^2}{t+k}\right].
\end{aligned} \tag{3.19}$$

Now, taking the sum over $k$ and $\rho_1$ and taking the expectation over the initial state, we find that the contribution to $E\{\sum_{j=1}^{R(t)-R(\gamma t)} |\mathcal{C}_j|^2\}$ coming from scenario (i) is at most

$$\frac{C_4 D}{2\sqrt{\gamma t}} \sum_{k\geq 0} \exp\biggl[-\frac{C_5 k^2}{t+k}\biggr] E\biggl\{\sum_{i\in[\sqrt{t},2\sqrt{t}]} N(i,0)\biggr\} \leq \frac{C_6\sqrt{t}}{\sqrt{\gamma}} \tag{3.20}$$

for a suitable constant $C_6$. By interchanging the roles of $\rho_1$ and $\rho_2$, we get the same contribution from scenario (ii).

The contribution from scenario (iii) can be estimated similarly. For fixed distinct $\rho_1 - \rho_3$ we now replace (3.18) by the bound

$$\begin{aligned} &\sum_{k\geq 0} \int_{\gamma t}^{t} P\{R^r(s-;\rho_1,\rho_2) = k, \\ &\qquad\qquad \pi(\rho_3, s-) = k+1 \text{ and } \rho_3 \text{ is white at time } s- | \mathcal{F}_0\} \\ &\qquad \times P\{\pi(\rho_1, s-) = \pi(\rho_2, s-) = k+1\} \frac{D}{2}\, ds \\ &\quad \leq \frac{C_4^2 D}{2\gamma t} \sum_{k\geq 0} P\{\text{in system without } \rho_1,\rho_2, R^r(\cdot;\rho_1,\rho_2) \text{ increases} \\ &\qquad\qquad\qquad \text{from } k \text{ to } k+1 \text{ by } 1 \text{ due to jump of } \rho_3 | \mathcal{F}_0\} \\ &\qquad\qquad \times \exp\biggl[-\frac{2C_5 k^2}{t+k}\biggr]. \end{aligned} \tag{3.21}$$

Analogously to (3.19), the sum of the right-hand side of (3.21) over $\rho_3$ is at most

$$\frac{C_4^2 D}{2\gamma t} \sum_{k\geq 0} \exp\biggl[-\frac{2C_5 k^2}{t+k}\biggr] \leq \frac{C_7}{\gamma\sqrt{t}}.$$

After summing over those $\rho_1, \rho_2$ which start in $[\sqrt{t}, 2\sqrt{t}]$ and taking expectation over the initial state, we find that scenario (iii) contributes at most

$$\frac{C_7}{\gamma\sqrt{t}} E\biggl\{\biggl[\sum_{i\in[\sqrt{t},2\sqrt{t}]} N(i,0)\biggr]^2\biggr\} \leq \frac{C_8}{\gamma}\sqrt{t}$$

to $E\{\sum_{j=1}^{R(t)-R(\gamma t)} |\mathcal{C}_j|^2\}$. Adding this to the contributions in (3.17) and (3.20), we find that

$$E\biggl\{\sum_{j=1}^{R(t)-R(\gamma t)} |\mathcal{C}_j|^2\biggr\} \leq \frac{C_9}{\gamma}\sqrt{t}. \tag{3.22}$$

Thus, (3.13) with $\delta = 4C_9/(\gamma\varepsilon)$ follows from Markov's inequality. $\square$

REMARK 4. Consider the system in which all particles perform asymmetric random walks, as described in Remark 2. Assume that $0 < p_+ < 1/2 < p_- < 1$ so that the particles have a drift to the left. It is intuitively clear that in this case, $R(t)$ should go to infinity at least linearly in $t$. In fact, by Remark 2 it cannot grow faster than linearly in $t$. It is possible to prove that there exists a constant $C_{10} > 0$ such that

$$\liminf \frac{1}{t} R(t) \geq C_{10} = C_{10}(p_+) \qquad \text{a.s.} \tag{3.23}$$

To prove (3.23), we use the same "second moment method" as in the proof of Theorem 3. To be precise, we use the following analogues of (3.12) and (3.22): this time let $\mathcal{C}_j$ be the collection of particles which turn black at the $j$th jump of $R(\cdot)$. Then for

$$\sum_{j=1}^{R(t)} |\mathcal{C}_j| \geq C_{11} t \tag{3.24}$$

and

$$\sum_{j=1}^{R(t)} |\mathcal{C}_j|^2 \leq C_{12} t. \tag{3.25}$$

If both these relations hold, then by Schwarz' inequality, as in the proof of Theorem 3,

$$C_{11}^2 t^2 \leq R(t) C_{12} t, \qquad \text{whence } R(t) \geq \frac{C_{11}^2}{C_{12}} t. \tag{3.26}$$

Now, (3.24) is trivial since all white particles have a drift $(p_- - p_+)$ toward the origin. Thus, the number of particles which reach the origin during $[0,t]$ is at least

$$\sum_{1 \leq x \leq (p_- - p_+)Dt/2} N(x,0) \geq \frac{\mu D}{4} (p_- - p_+) t,$$

outside an event of exponentially (in $t$) small probability. Again as in the proof of Theorem 3, all these particles will have coincided with the right edge of the aggregate and have changed color by time $t$, and are counted in $\sum_{j \leq R(t)} |\mathcal{C}_j|$. Thus, (3.24) with $C_{11} = \mu D(p_- - p_+)/4$ holds for all large $t$, a.s.

To prove (3.25), we need much of the machinery developed in Kesten and Sidoravicius (2003a) and we do not give this proof here.

**4. Positive recurrence in Caricature II.** In this section, we consider Caricature II, as described at the end of the Introduction. We denote the locations at time $t$ of the $J$ white particles relative to $R(t)$ as $X_1(t), X_2(t), \dots, X_J(t)$. We stress that these are the relative locations with respect to $R(t)$. The actual locations of the white particles in $\mathbb{Z}_+$ are $R(t)+X_1(t), \dots, R(t)+X_J(t)$. The process can be constructed in the following way. Let $\{S(t)\}_{t\geq 0}$ and $\{S_j(t)\}_{t\geq 0}, 1 \leq j \leq J$, be i.i.d. continuous-time simple random walks which start at $\mathbf{0}$ and have jump rate $D$. Also, let $\{Y_{j,k}, 1 \leq j \leq J, k \geq 1\}$ be i.i.d. random variables with common distribution $G$, concentrated on $\{1, 2, \dots\}$. The $\{Y_{j,k}\}$ are taken independent of the $\{S_j(t)\}$. Let $R(0) = \mathbf{0}, \tau_0 = 0$ and $X_j(0) = A_{j,0} \in \{1, 2, \dots\}$, $1 \leq j \leq J$. These $X_j(0)$ are regarded as nonrandom, but are otherwise arbitrary integers $\geq 1$. Then, when $\tau_k$ and $X_j(\tau_k)$ have already been determined for some $k \geq 0$, define

$$\tau_{k+1} = \inf\{t > \tau_k : X_j(\tau_k) + [S_j(t) - S_j(\tau_k)] = 0 \text{ for some } 1 \leq j \leq J\} \tag{4.1}$$

and let $r(k+1)$ be the value of $j$ for which $X_j(\tau_k) + [S_j(\tau_{k+1}) - S_j(\tau_k)] = 0$. Since almost surely only one of the random walks jumps at any time $t$, this index $r(k+1)$ is a.s. unique. $S_{r(k+1)}$ is the a.s. unique $S_j$ which has a jump at time $\tau_{k+1}$. Further, let

$$X_j(t) = X_j(\tau_k) + [S_j(t) - S_j(\tau_k)] \qquad \text{for } \tau_k \leq t < \tau_{k+1}, 1 \leq j \leq J \tag{4.2}$$

and

$$X_j(\tau_{k+1}) = X_j(\tau_{k+1}-) - I[j = r(k+1)] + A_{j,k+1}, \qquad 1 \leq j \leq J, \tag{4.3}$$

where the so-called adjustments $A_{j,k+1}$ are defined by

$$A_{j,k+1} = \begin{cases} Y_{r(k+1),k+1}, & \text{if } j = r(k+1), \\ Y_{j,k+1} - 1, & \text{if } X_j(\tau_{k+1}-) = 1, \text{ but } j \neq r(k+1), \\ -1, & \text{if } X_j(\tau_{k+1}-) \geq 2. \end{cases} \tag{4.4}$$

Note that a jump of a simple random walk is $+1$ or $-1$, so the $X_{r(k+1)}$, which jumps to 0 at time $\tau_{k+1}$, must satisfy

$$X_{r(k+1)}(\tau_{k+1}-) = 1. \tag{4.5}$$

Thus, the adjustments have been defined in such a way that

$$X_j(t) = S_j(t) + \sum_{k \geq 0 : \tau_k \leq t} A_{j,k}. \tag{4.6}$$

$R(t)$ is defined by (1.3). We point out that we make no adjustments at time 0, but start with the nonrandom $X_j(0) = A_{j,0}$, find $X_j(\tau_1)$ from (4.1)–(4.3) and then find $X_j(\tau_k)$ successively for $k = 2, 3, \dots$.

The preceding paragraph almost surely defines the $X_j$ for all time. It is clear from the description of the model that the (ordered) $J$-tuple $\{X_1(t), \dots, X_J(t)\}_{t\geq 0}$ is a strong Markov process. Its countable state space is $\{1, 2, \dots\}^J$. In a sequence of lemmas, we shall prove the following theorem.

THEOREM 4. *Consider Caricature II. If*

$$\mu_{10} := \sum_{n=1}^{\infty} n^{10} G(\{n\}) < \infty, \tag{4.7}$$

*then there exists a $J_0$ such that for all $J \geq J_0$, $\{X_j(t), 1 \leq j \leq J\}_{t \geq 0}$ is irreducible and positive recurrent. Moreover, under condition (4.7) and $J \geq J_0$, it holds almost surely that*

$$\lim_{t\to\infty} \frac{R(t)}{t} \text{ exists and is strictly greater than } 0. \tag{4.8}$$

We need more notation. $\ell_1(t), \dots, \ell_J(t)$ will be the values $X_1(t), \dots, X_J(t)$ in increasing order, so

$$\ell_1(t) \leq \ell_2(t) \leq \cdots \leq \ell_J(t). \tag{4.9}$$

We set $\widetilde{L}_0 = L_0 = \ell_1(0) + \cdots + \ell_J(0)$ and

$$\widetilde{L}_k := \ell_1(\tau_k-) + \cdots + \ell_J(\tau_k-) - 1, \qquad k \geq 1.$$

The values in the right-hand side here are the values of the $\ell_j$ "just before the $r(k)$th particle has jumped at time $\tau_k$ and before the adjustments at $\tau_k$ have been made." For later use, we note that

$$\ell_j(\tau_k-) \geq 1, \qquad 1 \leq j \leq J, \quad \text{and} \quad \widetilde{L}_k \geq J - 1 \tag{4.10}$$

because each $X_j(t) \geq 1$. We also define the $J$-vectors

$$U_k := (X_1(\tau_k-), \dots, X_J(\tau_k-)) \tag{4.11}$$

at this time. The values after the adjustments give us

$$L_k := \ell_1(\tau_k) + \cdots + \ell_J(\tau_k) = \widetilde{L}_k + \sum_{j=1}^{J} A_{j,k}.$$

We further define

$$\Lambda = \Lambda(\alpha) = \left\{(x_1, \dots, x_J) : x_j \in \{1, 2, \dots\}, \sum_{j=1}^{J} x_j \leq \alpha + 1\right\},$$

$$\nu_1 = \nu_1(\alpha) = \inf\{k \geq 1 : U_k \in \Lambda(\alpha)\} = \inf\{k \geq 1 : \widetilde{L}_k \leq \alpha\}$$

and

$$\nu_{n+1} = \nu_{n+1}(\alpha) = \inf\{k > \nu_n : U_k \in \Lambda(\alpha)\}.$$

We can now outline the proof of Theorem 4. The $J$-vectors $U_k, k \geq 1$, form a Markov chain with the countable state space

$$\Gamma := \left\{(x_1, x_2, \dots, x_J) : x_i \in \mathbb{Z}_+, \min_i x_i = 1\right\}.$$

The minimal coordinate of any $U_k = \min_i X_i(\tau_k-)$ must equal 1 because a particle can jump to 0 at time $\tau_k$ only if it is at 1 just before the jump. The Markov chain $\{U_k\}$ visits the *finite* set $\Lambda(\alpha)$ successively at the times $\nu_1, \nu_2, \ldots$. It is not hard to prove that the chain $\{U_k\}$ is irreducible (see below) and, in fact, even the embedded Markov chain $\{U_{\nu_i}\}_{i\geq 1}$ is irreducible. The latter has the finite set $\Gamma \cap \Lambda$ as state space and therefore has a unique invariant probability measure, $\rho$, say, on $\Gamma \cap \Lambda$. We shall prove that under condition (4.7), and for $J \geq J_0$ for a suitable $J_0 < \infty$, for any choice of the initial state $(X_1(0), \ldots, X_J(0))$, we have

$$\nu_2 < \infty \qquad \text{a.s. and } E\{\nu_2\} < \infty, \tag{4.12}$$

and

$$\tau_{\nu_2} < \infty \qquad \text{a.s. and } E\{\tau_{\nu_2}\} < \infty. \tag{4.13}$$

Let us write $E^\rho$ for the expectation when the Markov chain $\{U_{\nu_i}\}$ starts with the distribution $\rho$ for $U_{\nu_1}$. We shall prove that (4.12) and (4.13) also imply that

$$E^\rho\{\nu_2 - \nu_1\} < \infty \tag{4.14}$$

and

$$E^\rho\{\tau_{\nu_2} - \tau_{\nu_1}\} < \infty. \tag{4.15}$$

It then follows from the law of large numbers for Markov additive processes that

$$\begin{aligned} &1 \leq \lim_{n\to\infty} \frac{1}{n}\nu_n = E^\rho\{\nu_2 - \nu_1\} < \infty \quad \text{and} \\ &0 < \lim_{n\to\infty} \frac{1}{n}\tau_{\nu_n} = E^\rho\{\tau_{\nu_2} - \tau_{\nu_1}\} < \infty. \end{aligned} \tag{4.16}$$

A proof of (4.16) can easily be given by a slight generalization of Chung (1967), Theorems I.15.1 and 2. One must apply the argument there to the Markov chain $\{U_{\nu_i}\}$ with the deterministic function $f(\cdot)$ in Chung (1967) replaced by the random function $[\nu_{i+1} - \nu_i]$ for (4.14), or $\tau_{\nu_{i+1}} - \tau_{\nu_i}$ for (4.15).

It follows from (4.16) that

$$0 < \lim_{n\to\infty} \frac{\nu_n}{\tau_{\nu_{n+1}}} = \lim_{n\to\infty} \frac{\nu_{n+1}}{\tau_{\nu_n}} = \frac{E^\rho\{\nu_2 - \nu_1\}}{E^\rho\{\tau_{\nu_2} - \tau_{\nu_1}\}} < \infty.$$

But, for $\tau_{\nu_n} \leq t < \tau_{\nu_{n+1}}$, it holds that $\nu_n \leq R(t) < \nu_{n+1}$ [see (1.3)] and

$$\frac{\nu_n}{\tau_{\nu_{n+1}}} \leq \frac{R(t)}{t} \leq \frac{\nu_{n+1}}{\tau_{\nu_n}}$$

so that a.s.

$$0 < \lim_{t \to \infty} \frac{R(t)}{t} = \frac{E^\rho\{\nu_2 - \nu_1\}}{E^\rho\{\tau_{\nu_2} - \tau_{\nu_1}\}} < \infty. \tag{4.17}$$

This will prove (4.8).

To start on the details, let us take care of the irreducibility of $\{U_k\}$ and $\{U_{\nu_i}\}$, and the proof of (4.14) and (4.15) from (4.12) and (4.13). Let $x' = (x'_1, \dots, x'_J)$ and $y = (y_1, \dots, y_J)$ be points in $\Gamma$. Assume, without loss of generality, that $y_1 = 1$. Suppose $U_k = (x'_1, \dots, x'_J)$ at time $\tau_k-$ and that after the adjustments at time $\tau_k$, $(X_1(\tau_k), \dots, X_J(\tau_k)) = (x''_1, \dots, x''_J)$ with $x''_i \geq 1$, $1 \leq i \leq J$. It is clear that the random walks $S_i$, $1 \leq i \leq J$, can then move from $x''_i$ to $y_i$, for $1 \leq i \leq J$, in such a way that $x''_i + S_i$ stays $\geq 1$. Suppose that this happens over a time interval $[\tau_k, \tau_k + s)$ so that $X_i(\tau_k + s) = y_i, 1 \leq i \leq J$. Assume that the next jump of some $X_i$ occurs at time $\tau_k + s + u$ and that it is $X_1$ which jumps at that time from $X_1(\tau_k + s + u-) = y_1 = 1$ to 0. In this case, $\tau_{k+1} = \tau_k + s + u$ and $U_{k+1} = (y_1, \dots, y_J)$. Since $(y_1, \dots, y_J)$ is an arbitrary point in $\Gamma$, this proves the irreducibility of the chain $\{U_k\}$. In fact, it proves that $P\{U_{k+1} = y | U_k = x'\} > 0$ for any $x', y \in \Gamma$. We then automatically also have $P\{U_{\nu_{i+1}} = y | U_{\nu_i} = x\} > 0$ for any $x, y \in \Lambda \subset \Gamma$. Thus, $\{U_{\nu_i}\}$ is irreducible, as claimed.

As for (4.14) and (4.15), these are not immediately obvious because we have treated $X(\tau_0-)$ differently from the $X(\tau_k-)$ for $k \geq 1$ by not applying any adjustments at time $\tau_0 = 0$. However, the preceding paragraph shows that if we start in any (nonrandom) state $X(0) = x = (x_1, \dots, x_J)$, and $y = (y_1, \dots, y_J) \in \Gamma \cap \Lambda$, then $P\{\nu_1 = 2, U_2 = y | X(0) = x\} > 0$. We then also have

$$\begin{aligned} E\{\nu_2 | X(0) = x\} &\geq E\{\nu_2 - \nu_1 | X(0) = x\} \\ &\geq P\{\nu_1 = 2, U_2 = y | X(0) = x\} E\{\nu_2 - \nu_1 | U_{\nu_1} = y\}. \end{aligned}$$

Thus, (4.12) implies $E\{\nu_2 - \nu_1 | U_{\nu_1} = y\} < \infty$ for any $y \in \Gamma \cap \Lambda$ and this, in turn, implies (4.14), because there are only finitely many $y$ in $\Gamma \cap \Lambda$. In a similar way, one deduces (4.15) from (4.13).

We now start on the proof of (4.12) and (4.13). We define

$$\delta_k = \tau_{k+1} - \tau_k$$

and the $\sigma$-fields

$$\begin{aligned} \mathcal{F}(t) &= \sigma\text{-field generated by } \{S_j(s) : 1 \leq j \leq J, 0 \leq s \leq t\}, \\ \mathcal{G}_k &= \mathcal{F}(\tau_k) \vee \{Y_{j,n}, 1 \leq j \leq J, n \leq k-1\}. \end{aligned}$$

Note that the $Y_{j,k}$ are not included in the set of variables which generate $\mathcal{G}_k$. Thus, the information in $\mathcal{G}_k$ determines $\tau_k, X_j(\tau_k-)$ and $\widetilde{L}_k$, but not the adjustments $A_{j,k}$ or the values of $X_j(\tau_k), L_k$. The $\sigma$-field which also includes the information on $Y_{j,k}, 1 \le j \le J$, is

$$\mathcal{H}_k := \sigma\text{-field generated by } \mathcal{G}_k \vee \{Y_{j,k}, 1 \le j \le J\}.$$

Throughout, $D_i$ will denote various constants with values in $(0,\infty)$ which are independent of $J$. The same symbol $D_i$ may have different values in different formulae. For an event $A$, $I[A]$ denotes the indicator function of $A$. For real numbers $a, b$, we write $a \wedge b$ for $\min(a,b)$ and $a \vee b$ for $\max(a,b)$.

LEMMA 2. *Assume that $2q \in \{2,3,\dots\}$ and that $p \ge 0$ is such that*

$$\sum_{n=1}^{\infty} n^{2q+p} G(\{n\}) < \infty. \tag{4.18}$$

*Then, for all $\varepsilon > 0$, there exists a $J(q,\varepsilon)$ such that for all $J \ge J(q,\varepsilon)$ and $k \ge 0$, it holds that*

$$E\{\delta_k^q [X_j(\tau_k)]^p | \mathcal{G}_k\} \le \varepsilon \left(\frac{\widetilde{L}_k}{J}\right)^{2q-1} [X_j(\tau_k-)]^p. \tag{4.19}$$

*[For $k = 0$, we interpret $X_j(\tau_0-)$ as $X_j(0)$ and $J(q,\varepsilon)$ will depend also on $\ell_1$ if $k = 0$.]*

PROOF. We fix $k$ and abbreviate $\ell_j(\tau_k)$ to $\ell_j$. For the time being, we condition on $\mathcal{H}_k$ and consequently regard $\ell_j$, and also $\widetilde{L}_k, L_k$, as fixed. The main part of the proof is to show that

$$E\{\delta_k^q | \mathcal{H}_k\} \le \varepsilon \ell_1 \left(\frac{L_k}{J}\right)^{2q-1}. \tag{4.20}$$

The proof is based on the following well-known estimate: for $x \in \{1,2,\dots\}$,

$$P\{S_j(u) < -x\} \le \tfrac{1}{2} P\left\{\inf_{s \le u} S_j(s) \le -x\right\} \tag{4.21}$$

[see Doob (1953), proof of Theorem III.2.2]. This implies that

$$\begin{aligned} &P\{x_j + S_j(s) > 0 \text{ for } 0 \le s \le u\} \\ &\quad = P\left\{\inf_{s \le u} S_j(s) \ge 1 - x_j\right\} \\ &\quad = 1 - P\left\{\inf_{s \le u} S_j(s) \le -x_j\right\} \le 1 - 2P\{S_j(u) < -x_j\} \\ &\quad = P\{-x_j \le S_j(u) \le x_j\} \le 1 \wedge \frac{D_1 x_j}{\sqrt{u}}, \qquad u \ge 1, x_j \in \{1,2,\dots\} \end{aligned} \tag{4.22}$$

(by the local central limit theorem). Since $\delta_k > u$ occurs if and only if $X_j(t) = X_j(\tau_k) + [S_j(t) - S_j(\tau_k)] > 0$ for $\tau_k \le t \le \tau_k + u$, for $1 \le j \le J$, it follows that

$$P\{\delta_k > u | \mathcal{H}_k\} \le \prod_{j=1}^{J} \left[1 \wedge \frac{D_1 \ell_j}{\sqrt{u}}\right]. \tag{4.23}$$

For the remainder of this proof, we restrict ourselves to the case when $q$ is an integer. The case when $q$ is not an integer is actually easier. We set $\ell_0 = 0$ and $\ell_{J+1} = \infty$ and, without loss of generality, we assume $J \ge 8q + 2 \ge 10$. We further interpret the product $\ell_1 \cdots \ell_0$ as 1. Then (4.23) gives

$$\begin{aligned} E\{\delta_k^q | \mathcal{H}_k\} &= q \int_0^\infty u^{q-1} P\{\delta_k > u | \mathcal{H}_k\}\, du \\ &\le q \sum_{0 \le j < 2q} \ell_1 \cdots \ell_j [D_1]^j \int_{[D_1 \ell_j]^2}^{[D_1 \ell_{j+1}]^2} u^{q-1-j/2}\, du \\ &\quad + q \ell_1 \cdots \ell_{2q} [D_1]^{2q} \int_{[D_1 \ell_{2q}]^2}^{[D_1 \ell_{2q+1}]^2} u^{-1}\, du \\ &\quad + q \ell_1 \cdots \ell_{2q+1} [D_1]^{2q+1} \int_{[D_1 \ell_{2q+1}]^2}^{\infty} u^{-3/2}\, du \\ &\le D_2 \sum_{0 \le j < 2q} \ell_1 \cdots \ell_j [\ell_{j+1}]^{2q-j} \\ &\quad + D_2 \ell_1 \cdots \ell_{2q} \log \frac{\ell_{2q+1}}{\ell_{2q}} + D_2 \ell_1 \cdots \ell_{2q} \\ &\le D_3 \ell_1 \cdots \ell_{2q} \left[1 + \log \frac{\ell_{2q+1}}{\ell_{2q}}\right] \qquad (\text{since } \ell_1 \le \ell_2 \le \cdots \le \ell_J). \end{aligned} \tag{4.24}$$

It further follows from $\ell_1 \le \ell_2 \le \cdots \le \ell_J$ that

$$\ell_1 \le \ell_2 \le \cdots \le \ell_{2q+1} \le L_k/(J - 2q) \le 2L_k/J \qquad \text{for } J \ge 4q. \tag{4.25}$$

Substitution into (4.24) yields

$$\begin{aligned} E\{\delta_k^q | \mathcal{H}_k\} &\le D_3 \ell_1 \cdots \ell_{2q-1} \frac{L_k}{J} \left[\frac{J \ell_{2q}}{L_k} + \frac{J \ell_{2q}}{L_k} \log\left(\frac{2L_k}{J \ell_{2q}}\right)\right] \\ &\le D_3 \ell_1 [\ell_{2q}]^{2q-2} \frac{L_k}{J} \left[\frac{J \ell_{2q}}{L_k} + \frac{J \ell_{2q}}{L_k} \log\left(\frac{2L_k}{J \ell_{2q}}\right)\right]. \end{aligned} \tag{4.26}$$

Now, fix $0 < \eta = \eta(\varepsilon) < 2/e$ such that

$$D_3 \lambda^{2q-2} [\lambda + \lambda \log(2/\lambda)] \le \varepsilon \qquad \text{for } 0 \le \lambda \le \eta$$

and define the event

$$\mathcal{A}_k = \left\{\ell_{2q} \le \frac{\eta L_k}{J}\right\}.$$

Note that $\mathcal{A}_k \in \mathcal{H}_k$ and that the desired inequality (4.20) certainly holds on the event $\mathcal{A}_k$.

In order to deal with the conditional expectation on the complement of $\mathcal{A}_k$, we shall refine the estimate (4.24). First, we note that on $\mathcal{A}_k^c$, for all $j \le J/2$, it holds that

$$\ell_{2q} \ge \frac{\eta L_k}{J} \ge \frac{\eta}{J}(J-j+1)\ell_j \ge \frac{\eta}{2}\ell_j. \tag{4.27}$$

Thus, if we set $D_4 = (2/\eta) \vee 1$, then on $\mathcal{A}_k^c$,

$$\ell_j \le D_4 \ell_{2q}, \qquad 1 \le j \le J/2.$$

For such $\ell_j, \ell_{2q}$, there exists, for each $\zeta > 0$, some $f(\zeta) = f(\zeta, \varepsilon) > 0$ such that

$$\begin{aligned} P\{\ell_j + S_j(t) > 0 \text{ for } t \le [\zeta D_1 \ell_{2q}]^2\} &\le P\{S([\zeta D_1 \ell_{2q}]^2) > -D_4 \ell_{2q}\} \\ &\le 1 - f(\zeta), \qquad 1 \le j \le J/2. \end{aligned} \tag{4.28}$$

We shall soon fix a number $\zeta \in (0,1)$, but we need some inequalities before we can do so. If $0 < \zeta \le 1$ and $\ell_{2q-1}, \ell_{2q}$ are such that $\zeta \ell_{2q} \ge \ell_{2q-1}$, then we replace (4.24) on $\mathcal{A}_k^c$ by (recall $\ell_0 = 0$)

$$\begin{aligned} E\{\delta_k^q | \mathcal{H}_k\} &\le q \sum_{0 \le j \le 2q-2} D_1^j \ell_1 \cdots \ell_j [[D_1 \ell_{j+1}]^{2q-j} - [D_1 \ell_j]^{2q-j}] \\ &\quad + D_2 \ell_1 \cdots \ell_{2q-1} \int_{[D_1 \ell_{2q-1}]^2}^{[\zeta D_1 \ell_{2q}]^2} u^{q-1-(2q-1)/2} \, du \\ &\quad + D_2 \ell_1 \cdots \ell_{2q-1} \\ &\quad\quad \times \int_{[\zeta D_1 \ell_{2q}]^2}^{[D_1 D_4 \ell_{2q}]^2} u^{q-1-(2q-1)/2} [1 - f(\zeta)]^{\lfloor J/2 \rfloor - 2q+1} \, du \\ &\quad + D_2 \ell_1 \cdots \ell_{2q-1} \\ &\quad\quad \times \int_{[D_1 D_4 \ell_{2q}]^2}^{\infty} u^{q-1-q+1/2} \left[\frac{D_1 D_4 \ell_{2q}}{\sqrt{u}}\right]^{\lfloor J/2 \rfloor - 2q+1} du \\ &\le D_5 \ell_1 \cdots \ell_{2q-2} [\ell_{2q-1}]^2 + D_5 \zeta \ell_1 \cdots \ell_{2q} \\ &\quad + 2 D_1 D_2 D_4 \ell_1 \cdots \ell_{2q} [1 - f(\zeta)]^{\lfloor J/2 \rfloor - 2q+1} \\ &\quad + D_1 D_2 D_4 \ell_1 \cdots \ell_{2q} [\lfloor J/2 \rfloor / 2 - q]^{-1} \\ &\le D_6 \ell_1 \cdots \ell_{2q} [\zeta + [1 - f(\zeta)]^{\lfloor J/2 \rfloor - 2q+1} + J^{-1}]. \end{aligned} \tag{4.29}$$

Here, we used the fact that $J \ge 8q + 2$. We also used (4.28) to conclude that

$$P\{\ell_j + S_j(s) > 0 \text{ for } s \le u\} \le \min\left\{\frac{D_1 \ell_j}{\sqrt{u}}, 1 - f(\zeta)\right\}$$

for all $j \le J/2$ and $u$ in the domain of integration in the second and third integral in the first right-hand side (4.29). The constant $D_6$ depends on $q$ and $\eta$ (or $\varepsilon$) only and not on $\zeta$ or $J$, provided $J \ge 8q+2$. Without loss of generality, we take $\varepsilon \le 1$ and $D_6 \ge 1$.

Finally, we take

$$\zeta = \frac{\varepsilon}{3 \cdot 2^{2q-1} D_6} \le 1$$

and use (4.25). We then see from (4.29) that there exists an $J = J(q,\varepsilon)$ such that for $J \ge J(q,\varepsilon)$ on $\mathcal{A}^c \cap \{\ell_{2q-1} \le \zeta \ell_{2q}\}$,

$$E\{\delta_k^q | \mathcal{H}_k\} \le D_6 \ell_1 \left(\frac{2L_k}{J}\right)^{2q-1} 3\zeta \le \varepsilon \ell_1 \left(\frac{L_k}{J}\right)^{2q-1}.$$

If $\ell_{2q-1} > \zeta \ell_{2q}$, then let

$$j_0 = \max\{j \ge 0 : \ell_j \le \zeta \ell_{2q}\}.$$

Recall that $\ell_0 = 0$, so $j_0$ is well defined. Also, $j_0 \le 2q-2$ since $\ell_{2q} \ge \ell_{2q-1} > \zeta \ell_{2q}$. Instead of (4.29), we now use

$$\begin{aligned}
E\{\delta_k^q | \mathcal{H}_k\} \le{}& q \sum_{0 \le j \le j_0 - 1} D_1^j \ell_1 \cdots \ell_j [[D_1 \ell_{j+1}]^{2q-j} - [D_1 \ell_j]^{2q-j}] \\
&+ D_2 \ell_1 \cdots \ell_{j_0} \int_{[D_1 \ell_{j_0}]^2}^{[\zeta D_1 \ell_{2q}]^2} u^{q-1-j_0/2}\, du \\
&+ D_2 \ell_1 \cdots \ell_{2q-2} \\
&\times \int_{[\zeta D_1 \ell_{2q}]^2}^{[D_1 D_4 \ell_{2q}]^2} u^{q-1-(2q-2)/2} [1 - f(\zeta)]^{\lfloor J/2 \rfloor - 2q + 2}\, du \\
&+ D_2 \ell_1 \cdots \ell_{2q-1} \\
&\times \int_{[D_1 D_4 \ell_{2q}]^2}^{\infty} u^{q-1-q+1/2} \left[\frac{D_1 D_4 \ell_{2q}}{\sqrt{u}}\right]^{\lfloor J/2 \rfloor - 2q + 1} du \\
\le{}& D_5 \zeta \ell_1 [\ell_{2q}]^{2q-1} \\
&+ D_5 \ell_1 \cdots \ell_{2q-2} [\ell_{2q}]^2 [1 - f(\zeta)]^{\lfloor J/2 \rfloor - 2q + 2} \\
&+ D_1 D_2 D_4 \ell_1 \cdots \ell_{2q} [\lfloor J/2 \rfloor / 2 - q]^{-1} \\
\le{}& D_6 \ell_1 [\ell_{2q}]^{2q-1} [\zeta + \zeta^{-1} [1 - f(\zeta)]^{\lfloor J/2 \rfloor - 2q + 2} + J^{-1}].
\end{aligned} \tag{4.30}$$

The factor $\zeta^{-1}$ multiplying $[1 - f(\zeta)]^{\lfloor J/2 \rfloor - 2q + 2}$ in the last member of (4.30) has been inserted to deal with the case $q = 1$. Note that (4.30) is also valid in the case $j_0 = 0$ (which contains the case $q = 1, \zeta \ell_2 < \ell_1$). Indeed, $j_0 = 0$ means that $\zeta\ell_{2q} < \ell_1$ and then the sum in the first right-hand side of (4.30) is empty, while the first integral becomes

$$D_2 \int_0^{[\zeta D_1 \ell_{2q}]^2} u^{q-1}\,du = \frac{D_2}{q}[\zeta D_1 \ell_{2q}]^{2q} \le \frac{D_2}{q} D_1 \ell_1 [\zeta D_1 \ell_{2q}]^{2q-1}.$$

We leave it to the reader to check that in all these cases, there exists a $J(q,\varepsilon)$ such that for $J \ge J(q,\varepsilon)$, on $\mathcal{A}^c \cap \{\ell_{2q-1} > \zeta\ell_{2q}\}$, (4.20) holds so that (4.20) has been proven in general.

To go from (4.20) to (4.19), we take conditional expectations with respect to $\mathcal{G}_k$, which is a sub-$\sigma$-field of $\mathcal{H}_k$. This gives

$$\begin{aligned} E\{\delta_k^q [X_j(\tau_k)]^p | \mathcal{G}_k\} &= E\{E\{\delta_k^q [X_j(\tau_k)]^p | \mathcal{H}_k\} | \mathcal{G}_k\} \\ &\le \varepsilon E\left\{\ell_1 \left(\frac{L_k}{J}\right)^{2q-1} [X_j(\tau_k)]^p | \mathcal{G}_k\right\}. \end{aligned} \tag{4.31}$$

Taking this conditional expectation amounts to integrating out the $Y_{j,k}$, which are independent of $\mathcal{G}_k$. For the remainder of this proof, we restrict ourselves to the case $k \ge 1$; when $k = 0$, the proof simplifies. We use the facts that

$$0 \le X_j(\tau_k) = X_j(\tau_k-) - I[j = r(k)] + A_{j,k} \le X_j(\tau_k-) + Y_{j,k} \tag{4.32}$$

[see (4.3)] and

$$L_k = \widetilde{L}_k + \sum_{j=1}^J A_{j,k} \le \widetilde{L}_k + \sum_{j=1}^J Y_{j,k}, \tag{4.33}$$

and $\ell_1(\tau_k) \le Y_{r(k),k}$ since $\ell_1(\tau_k)$ is the minimum of the $X_j(\tau_k), 1 \le j \le J$, and $X_{r(k)}(\tau_k) = Y_{r(k),k}$ by (4.3)–(4.5). Note that both $r(k)$ and $\widetilde{L}_k$ are $\mathcal{G}_k$-measurable. Moreover, $Y_{r(k),k}$ is independent of $\mathcal{G}_k$ (and hence of $\widetilde{L}_k$) and has the distribution $G$. Also, for $j \ne r(k)$, $Y_{j,k}$ and $Y_{r(k),k}$ are independent. Thus,

$$\begin{aligned} &E\{\ell_1(L_k)^{2q-1} [X_j(\tau_k)]^p | \mathcal{G}_k\} \\ &\quad \le E\left\{Y_{r(k),k} \left[\widetilde{L}_k + \sum_{\ell=1}^J Y_{\ell,k}\right]^{2q-1} [X_j(\tau_k-) + Y_{j,k}]^p | \mathcal{G}_k\right\}. \end{aligned} \tag{4.34}$$

We shall frequently use the following special case of Hölder's inequality: for any $a_j \ge 0$ and for $p \ge 1$,

$$\left[\sum_{j=1}^n a_j\right]^p \le n^{p-1} \sum_{j=1}^n a_j^p. \tag{4.35}$$

In particular, as a case with $n=2$ and with $n=J$, we have

$$\begin{aligned}\left[\widetilde{L}_k + \sum_{\ell=1}^J Y_{\ell,k}\right]^{2q-1} &\le 2^{2q-2}\left\{[\widetilde{L}_k]^{2q-1} + \left[\sum_{\ell=1}^J Y_{\ell,k}\right]^{2q-1}\right\}\\ &\le 2^{2q-2}[\widetilde{L}_k]^{2q-1} + (2J)^{2q-2}\sum_{\ell=1}^J [Y_{\ell,k}]^{2q-1}.\end{aligned}$$

Also, for $p \ge 0$,

$$[X_j(\tau_k-) + Y_{j,k}]^p \le 2^p[X_j(\tau_k-)]^p + 2^p[Y_{j,k}]^p. \tag{4.36}$$

(We have a factor $2^p$ instead of $2^{p-1}$ here to deal with the case $0 \le p < 1$.) In agreement with the notation of (4.7), we write $\mu_\kappa$ for the $\kappa$th moment of $G$ and use the fact that $E\{[Y_{\ell,k}]^a[Y_{m,k}]^b\} \le \mu_{a+b}$ for any $1 \le \ell, m \le J$ and $a, b \ge 0$, by Hölder's inequality. We can therefore continue (4.31) with

$$\begin{aligned}&E\{\ell_1(L_k)^{2q-1}[X_j(\tau_k)]^p|\mathcal{G}_k\}\\ &\quad\le 2^{2q+p-2}[\widetilde{L}_k]^{2q-1}[X_j(\tau_k-)]^p E\{Y_{r(k),k}\}\\ &\qquad + 2^{2q+p-2}[\widetilde{L}_k]^{2q-1}E\{Y_{r(k),k}[Y_{j,k}]^p\}\\ &\qquad + 2^{2q+p-2}J^{2q-2}[X_j(\tau_k-)]^p\sum_{\ell=1}^J E\{Y_{r(k),k}[Y_{\ell,k}]^{2q-1}\}\\ &\qquad + 2^{2q+p-2}J^{2q-2}\sum_{\ell=1}^J E\{Y_{r(k),k}[Y_{\ell,k}]^{2q-1}[Y_{j,k}]^p\}\\ &\quad\le 2^{2q+p-2}\mu_1[\widetilde{L}_k]^{2q-1}[X_j(\tau_k-)]^p + 2^{2q+p-2}\mu_{p+1}[\widetilde{L}_k]^{2q-1}\\ &\qquad + 2^{2q+p-2}J^{2q-1}[X_j(\tau_k-)]^p\mu_{2q} + 2^{2q+p-2}J^{2q-1}\mu_{2q+p}.\end{aligned}$$

Substitution of this estimate into (4.31) shows that if we take $J \ge J(q,\varepsilon) \ge 2$, then for $J \ge J(q,\varepsilon)$,

$$\begin{aligned}E\{\delta_k^q[X_j(\tau_k)]^p|\mathcal{G}_k\} &\le \varepsilon 2^{2q+p-2}\left(\frac{\widetilde{L}_k}{J}\right)^{2q-1}[\mu_1[X_j(\tau_k-)]^p + \mu_{p+1}]\\ &\qquad + \varepsilon 2^{2q+p-2}[[X_j(\tau_k-)]^p\mu_{2q} + \mu_{2q+p}]\\ &\le D_7\varepsilon\left(\frac{\widetilde{L}_k}{J}\right)^{2q-1}[X_j(\tau_k-)]^p\end{aligned} \tag{4.37}$$

[recall that $X_j(\tau_k-) \ge 1$ and $\widetilde{L}_k \ge J-1$; see (4.10)]. The lemma follows by replacing $\varepsilon$ by $\varepsilon/D_7$. $\square$

We define

$$Q_q(t) = \sum_{j=1}^{J} [X_j(t)]^q, \qquad \widetilde{Q}_{q,0} = \sum_{j=1}^{J} [X_j(0)]^q \quad \text{and for } k \geq 1,$$

$$\widetilde{Q}_{q,k} = \sum_{j=1}^{J} [X_j(\tau_k -) - I[j = r(k)]]^q.$$

The quantity $\widetilde{Q}_{q,k}$ is a $q$th power analog of $\widetilde{L}_k$. $\widetilde{Q}_{q,k}$ is the sum of the $q$th powers of the coordinates "just before" the adjustments at $\tau_k$, but taking into account the jump of one coordinate from 1 to 0 at $\tau_k$.

Lemma 3. *Let the $X_j(0) \geq 1$ be fixed (nonrandom) and let $q \in \{2, 3, \dots\}$ and $k \geq 1$. Assume that $\mu_q := \sum_{n=1}^{\infty} n^q G(\{n\}) < \infty$. There then exists some $J(q)$ and for $J \geq J(q)$, there exists an $\alpha(q) = \alpha(q, J)$ such that for $J \geq J(q)$ and $\alpha \geq \alpha(q, J)$,*

$$\begin{aligned} E\{\widetilde{Q}_{q,k+1} | \mathcal{G}_k\} &\leq \widetilde{Q}_{q,k} - q 2^{-q-1} \widetilde{Q}_{(q-1),k} \\ &\leq \widetilde{Q}_{q,k} - 1 \qquad \text{on the event } \{\widetilde{L}_k > \alpha\}. \end{aligned} \tag{4.38}$$

*Consequently, $\nu_1(\alpha) < \infty$ a.s. and the process $\{X(\tau_k -)\}_{k \geq 0}$ is recurrent. Also, there exists some constant $C_1 = C_1(J, q, \alpha, X(0)) < \infty$ such that for $J \geq J(q)$ and $k \geq 0$,*

$$E\{\widetilde{Q}_{q,k+1} I[\nu_2(\alpha) > k]\} \leq C_1. \tag{4.39}$$

Proof. First, observe that all $\tau_\ell$ are almost surely finite. This follows by induction on $\ell$ from the fact that

$$\tau_\ell = \sum_{k=1}^{\ell-1} \delta_k \tag{4.40}$$

and the estimate (4.23) for the tail of the conditional distribution of $\delta_k$. If $\tau_\ell < \infty$, then almost surely $\ell_j(\tau_\ell) < \infty$, and then $\delta_\ell < \infty$ by (4.23), and hence $\tau_{\ell+1} < \infty$.

We shall need the following inequality. For $v \geq 2$, there exists a constant $D_8 = D_8(v, D)$ such that on the event $\{\tau_k < \infty\}$,

$$E\{|S_j(\tau_{k+1}) - S_j(\tau_k)|^v | \mathcal{H}_k\} \leq D_8[1 + E\{\delta_k^{v/2} | \mathcal{H}_k\}]. \tag{4.41}$$

This inequality is probably well known. For completeness, we shall outline a proof. Introduce

$$Z(s, n) = S_j\left(\left(\tau_k + \frac{s+1}{n}\right) \wedge \tau_{k+1}\right) - S_j\left(\left(\tau_k + \frac{s}{n}\right) \wedge \tau_{k+1}\right).$$

Then

$$S_j(\tau_{k+1}) - S_j(\tau_k) = \lim_{N\to\infty} \lim_{n\to\infty} \sum_{s=0}^{Nn} Z(s,n). \tag{4.42}$$

Moreover, for fixed $n$, the $Z(s,n), s \geq 0$, form a sequence of martingale differences with respect to the $\sigma$-fields

$$\mathcal{F}\biggl(\biggl(\tau_k + \frac{s}{n}\biggr) \wedge \tau_{k+1}\biggr) \vee \mathcal{H}_k.$$

Thus, by Fatou's lemma and the Burkholder–Davis inequality [see Gut (1988), Theorem A.2.2, or Hall and Heyde (1980), Theorem 2.10],

$$\begin{aligned}
E\{|S_j(\tau_{k+1}) - S_j(\tau_k)|^v | \mathcal{H}_k\} &\leq \liminf_{N\to\infty} \liminf_{n\to\infty} E\Biggl\{\Biggl|\sum_{s=0}^{Nn} Z(s,n)\Biggr|^v \Bigg| \mathcal{H}_k\Biggr\} \\
&\leq D_8(v) \liminf_{N\to\infty} \liminf_{n\to\infty} E\Biggl\{\Biggl|\sum_{s=0}^{Nn} Z^2(s,n)\Biggr|^{v/2} \Bigg| \mathcal{H}_k\Biggr\}.
\end{aligned}$$

But,

$$|Z(s,n)| \leq \text{number of jumps of } S_j \text{ during } ((\tau_k + s/n), (\tau_k + (s+1)/n)],$$

so

$$\begin{aligned}
\sum_{s=0}^{Nn} Z^2(s,n) &\leq \Biggl[\sum_{s=0}^{Nn} |Z(s,n)|\Biggr]^2 \\
&\leq [\text{number of jumps of } S_j \text{ during } (\tau_k, (\tau_k + N + 1))]^2.
\end{aligned} \tag{4.43}$$

Moreover,

$$\begin{aligned}
\lim_{n\to\infty} \sum_{s=0}^{Nn} Z^2(s,n) &= \text{number of jumps of } S_j \\
&\qquad \text{during } (\tau_k, (\tau_k + N) \wedge \tau_{k+1}].
\end{aligned} \tag{4.44}$$

Conditionally on $\mathcal{H}_k$, the number of jumps in the right-hand side of (4.43) is a Poisson variable with mean $D(N+1)$. Since a Poisson variable has all moments, it is not hard to show the equality

$$\begin{aligned}
&\liminf_{n\to\infty} E\Biggl\{\Biggl|\sum_{s=0}^{Nn} Z^2(s,n)\Biggr|^{v/2} \Bigg| \mathcal{H}_k\Biggr\} \\
&\qquad = E\Biggl\{\lim_{n\to\infty} \Biggl|\sum_{s=0}^{Nn} Z^2(s,n)\Biggr|^{v/2} \Bigg| \mathcal{H}_k\Biggr\}
\end{aligned}$$

$$
\begin{aligned}
&= E\{[\text{number of jumps of } S_j \text{ during} \\
&\qquad (\tau_k, (\tau_k + N) \wedge \tau_{k+1}]]^{v/2} | \mathcal{H}_k\} \\
&\le E\{[\text{number of jumps of } S_j \text{ during} \\
&\qquad (\tau_k, \tau_k + N \wedge (\lceil \tau_{k+1} - \tau_k \rceil)]]^{v/2} | \mathcal{H}_k\}.
\end{aligned} \tag{4.45}
$$

But, conditionally on $\mathcal{H}_k$, the jumps of $S_j$ during $(\tau_k, \infty)$ form a Poisson process with jump rate $D$. By writing the number of jumps in $(\tau_k, \tau_k + N \wedge (\lceil \tau_{k+1} - \tau_k \rceil)]$ as the sum over $r$ from 1 to $N \wedge (\lceil \tau_{k+1} - \tau_k \rceil)$ of the number of jumps in $(r-1, r]$ and using Theorem I.5.2 of Gut (1988), we then see that the right-hand side of (4.45) is at most

$$
D_9 E\{[N \wedge \lceil \tau_{k+1} - \tau_k \rceil]^{v/2} | \mathcal{H}_k\} \le D_9 2^{v-1} [1 + E\{[\tau_{k+1} - \tau_k]^{v/2} | \mathcal{H}_k\}]. \tag{4.46}
$$

The inequality (4.41) follows from (4.42)–(4.46) because $\delta_k = \tau_{k+1} - \tau_k$.

We now fix $q$ and $\alpha$. Before we start on the proof proper of (4.38), we should show that the conditional expectation in the left-hand side of (4.38) makes sense, that is, $E\{\widetilde{Q}_{q,k+1}\} < \infty$. To this end, we observe that by (4.2), on the event $\{\tau_k < \infty\}$,

$$
X_j(\tau_{k+1}-) - I[j = r(k+1)] = X_j(\tau_k) + S_j(\tau_{k+1}) - S_j(\tau_k)
$$

and hence

$$
\begin{aligned}
\widetilde{Q}_{q,k+1} &= \sum_{j=1}^{J} [X_j(\tau_k)]^q \\
&\quad + \sum_{j=1}^{J} \sum_{u=1}^{q} \binom{q}{u} [X_j(\tau_k)]^{q-u} [S_j(\tau_{k+1}) - S_j(\tau_k)]^u.
\end{aligned} \tag{4.47}
$$

Furthermore, for $k \ge 1$,

$$
X_j(\tau_k) = X_j(\tau_k-) - I[j = r(k)] + A_{j,k} \tag{4.48}
$$

[see (4.3)], so for $k \ge 1$ [see (4.35)],

$$
\begin{aligned}
\sum_{j=1}^{J} [X_j(\tau_k)]^q &= \sum_{j=1}^{J} [X_j(\tau_k-) - I[j = r(k)] + A_{j,k}]^q \\
&\le 2^{q-1} \sum_{j=1}^{J} [X_j(\tau_k-) - I[j = r(k)]]^q + 2^{q-1} \sum_{j=1}^{J} |A_{j,k}|^q \\
&= 2^{q-1} \widetilde{Q}_{q,k} + 2^{q-1} \sum_{j=1}^{J} |A_{j,k}|^q.
\end{aligned} \tag{4.49}
$$

Now, note that by (4.41),

$$E\{|S_j(\tau_1) - S_j(\tau_0)|^q\} \le D_8[1 + E\{\delta_1^{q/2}\}]$$

and by virtue of (the proof of) Lemma 2, there exists some $J_1(q)$ such that $E\{\delta_1^{q/2}\} < \infty$ for $J \ge J_1$. In fact, (4.23) and (4.24) show that $J > 2q$ suffices for this. Then (4.47) shows that $E\{\widetilde{Q}_{q,1}\} < \infty$ for $J \ge J_1$ [recall that the $X_j(0)$ are nonrandom]. From there on, we apply induction on $k$ to show that $E\{\widetilde{Q}_{q,k}\} < \infty$ and $E\{\delta_k^{q/2}\} < \infty$ for all $k$, by means of (4.47), (4.49), (4.41), (4.19) and

$$\begin{aligned} &E\{[X_j(\tau_k)]^{q-u}[S_j(\tau_{k+1}) - S_j(\tau_k)]^u\} \\ &\quad \le [E[X_j(\tau_k)]^q]^{(q-u)/q}[E\{|S_j(\tau_{k+1}) - S_j(\tau_k)|^q\}]^{u/q}. \end{aligned} \tag{4.50}$$

This shows that the conditional expectation in (4.38) is well defined.

We turn to the proof of (4.38) itself. The basic relation is

$$\begin{aligned} &E\{X^j(\tau_{k+1}-) - I[j = r(k+1)] - X^j(\tau_k)|\mathcal{H}_k\} \\ &\quad = E\{S^j(\tau_{k+1}) - S^j(\tau_k)|\mathcal{H}_k\} = 0, \end{aligned} \tag{4.51}$$

which follows from Wald's equation [see Chow and Teicher (1978), Theorem 5.3.1; this reference deals with discrete-time random walks only, but we can again approximate $\tau_{k+1}$ by $\lceil m\tau_{k+1}\rceil/m$ and let $m$ go to infinity]. Combined with (4.47), (4.51) implies that

$$\begin{aligned} E\{\widetilde{Q}_{q,k+1}|\mathcal{H}_k\} = \sum_{j=1}^J E\Bigg\{&[X_j(\tau_k)]^q + \sum_{u=2}^q \binom{q}{u} [X_j(\tau_k)]^{q-u} \\ &\times [S_j(\tau_{k+1}) - S_j(\tau_k)]^u |\mathcal{H}_k\Bigg\}. \end{aligned} \tag{4.52}$$

We shall fix $\varepsilon > 0$ in (4.60) below. In the last sum, we then have for $u \ge 2$ and $J \ge J(q, \varepsilon)$,

$$\begin{aligned} &\sum_{j=1}^J E\{[X_j(\tau_k)]^{q-u}|S_j(\tau_{k+1}) - S_j(\tau_k)|^u|\mathcal{H}_k\} \\ &\quad \le D_8 \sum_{j=1}^J [X_j(\tau_k)]^{q-u}[1 + E\{\delta_k^{u/2}|\mathcal{H}_k\}] \end{aligned} \tag{4.53}$$

[by (4.41) and the $\mathcal{H}_k$-measurability of $X_j(\tau_k)$; note that we used $X_j(t) \ge 0$ so that $X_j(t) = |X_j(t)|$]. Next, by (4.32) and (4.36), for $q \ge u$,

$$[X_j(\tau_k)]^{q-u} \le 2^{q-u}[X_j(\tau_k-)]^{q-u} + 2^{q-u}[Y_{j,k}]^{q-u}.$$

By taking conditional expectation with respect to $\mathcal{G}_k$ in (4.53) and using (4.19), we now obtain

$$\sum_{j=1}^{J} E\{[X_j(\tau_k)]^{q-u}|S_j(\tau_{k+1}) - S_j(\tau_k)|^u|\mathcal{G}_k\}$$
$$\le D_{10}\left[1+\varepsilon\left(\frac{\widetilde{L}_k}{J}\right)^{u-1}\right]\left[\sum_{j=1}^{J}[X_j(\tau_k-)]^{q-u} + D_{11}J\right].$$

But, $X_j(\tau_k-) = X_j(\tau_k-) - I[j = r(k)]$ for all $j \ne r(k)$, and $X_{r(k)}(\tau_k-) = 1$. Therefore,

$$J - 1 \le \sum_{j=1}^{J}[X_j(\tau_k-)]^{q-u} \le \widetilde{Q}_{(q-u),k} + 1 \le 2\widetilde{Q}_{(q-u),k}.$$

We then also have

$$\begin{aligned}
&E\{\widetilde{Q}_{q,k+1}|\mathcal{G}_k\}\\
&\quad = E\{\widetilde{Q}_{q,k+1} - Q_q(\tau_k)|\mathcal{G}_k\} + E\{Q_q(\tau_k) - \widetilde{Q}_{q,k}|\mathcal{G}_k\} + \widetilde{Q}_{q,k}\\
&\quad \le D_{12}\sum_{u=2}^{q}\left[1+\varepsilon\left(\frac{\widetilde{L}_k}{J}\right)^{u-1}\right]\widetilde{Q}_{(q-u),k}\\
&\qquad + E\{Q_q(\tau_k) - \widetilde{Q}_{q,k}|\mathcal{G}_k\} + \widetilde{Q}_{q,k}\\
&\quad = D_{12}\sum_{u=2}^{q}\left[1+\varepsilon\left(\frac{\widetilde{L}_k}{J}\right)^{u-1}\right]\widetilde{Q}_{(q-u),k}\\
&\qquad + \sum_{j=1}^{J} E\{[X_j(\tau_k-) - I[j = r(k)] + A_{j,k}]^q\\
&\qquad\qquad - [X_j(\tau_k-) - I[j = r(k)]]^q|\mathcal{G}_k\} + \widetilde{Q}_{q,k}.
\end{aligned} \tag{4.54}$$

We turn to a bound when $k \ge 1$ for the second sum in the right-hand side here. Its summand equals

$$\begin{aligned}
&E\{[X_j(\tau_k-) - I[j = r(k)] + A_{j,k}]^q - [X_j(\tau_k-) - I[j = r(k)]]^q|\mathcal{G}_k\}\\
&\quad = qE\{A_{j,k}[X^*(j,k)]^{q-1}|\mathcal{G}_k\}
\end{aligned} \tag{4.55}$$

for some $X^*(j,k)$ between $X_j(\tau_k-) - I[j = r(k)]$ and $X_j(\tau_k)$. We now consider two cases.

*Case* (i): $X_j(\tau_k-) \ge 2$. In this case, $A_{j,k} = -1$, $j \ne r(k)$, and

$$\begin{aligned}
\tfrac{1}{2}[X_j(\tau_k-) - I[j = r(k)]] &= \tfrac{1}{2}X_j(\tau_k-) \le X_j(\tau_k-) - I[j = r(k)] - 1\\
&= X_j(\tau_k) < X^*(j,k) < X_j(\tau_k-) - I[j = r(k)],
\end{aligned}$$

so

$$A_{j,k}[X^*(j,k)]^{q-1} \le -2^{1-q}[X_j(\tau_k-) - I[j = r(k)]]^{q-1}.$$

Hence, the right-hand side of (4.55) is at most

$$-q2^{1-q}[X^j(\tau_k-) - I[j = r(k)]]^{q-1}$$

in case (i).

*Case* (ii): $X_j(\tau_k-) = 1$. Then $0 \le A_{j,k} \le Y_{j,k}$ and

$$0 \le X_j(\tau_k-) - I[j = r(k)] \le X^*(j,k) \le X_j(\tau_k) = Y_{j,k}.$$

Thus, in this case,

$$A_{j,k}[X^*(j,k)]^{q-1} \le [Y_{j,k}]^q.$$

If we set

$$\widetilde{\mathcal{U}}_k = \{1 \le j \le J : X_j(\tau_k-) \ge 2\},$$

then

$$\begin{aligned}
&\sum_{j=1}^{J} E\{[X_j(\tau_k-) - I[j = r(k)] + A_{j,k}]^q - [X_j(\tau_k-) - I[j = r(k)]]^q | \mathcal{G}_k\} \\
&\quad \le q \sum_{j \in \widetilde{\mathcal{U}}_k} (-2^{1-q})[X_j(\tau_k-) - I[j = r(k)]]^{q-1} + q \sum_{j \notin \widetilde{\mathcal{U}}_k} E\{[Y_{j,k}]^q\} \\
&\quad \le (-q2^{1-q})\widetilde{Q}_{(q-1),k} + q2^{1-q} \sum_{j \notin \widetilde{\mathcal{U}}_k} [X_j(\tau_k-) - I[j = r(k)]]^{q-1} + qJ\mu_q \\
&\quad \le (-q2^{1-q})\widetilde{Q}_{(q-1),k} + q2^{1-q}(J - |\widetilde{\mathcal{U}}_k|) + qJ\mu_q \\
&\qquad\qquad\qquad (\text{because } X_j(\tau_k-) - I[j = r(k)] \le 1 \text{ for } j \notin \widetilde{\mathcal{U}}_k) \\
&\quad \le (-q2^{1-q})\widetilde{Q}_{(q-1),k} + D_{13}J.
\end{aligned}$$

Substitution of the estimates in cases (i) and (ii) into (4.54) shows that

$$\begin{aligned}
E\{\widetilde{Q}_{q,k+1} | \mathcal{G}_k\} &\le D_{12} \sum_{u=2}^{q} \left[1 + \varepsilon \left(\frac{\widetilde{L}_k}{J}\right)^{u-1}\right] \widetilde{Q}_{(q-u),k} \\
&\quad - q2^{1-q}\widetilde{Q}_{(q-1),k} + D_{13}J + \widetilde{Q}_{q,k}.
\end{aligned} \tag{4.56}$$

On $\{\widetilde{L}_k > \alpha\}$, if $q \ge 2$, we have

$$\alpha^{q-1} < [\widetilde{L}_k]^{q-1} \le J^{q-2}\widetilde{Q}_{(q-1),k}, \tag{4.57}$$

by virtue of (4.35). Thus, we can choose $\alpha > J$ so large that $D_{13}J \le q2^{-q}\widetilde{Q}_{(q-1),k}$ and

$$E\{\widetilde{Q}_{q,k+1}|\mathcal{G}_k\} \le D_{12}\sum_{u=2}^{q}\left[1+\varepsilon\left(\frac{\widetilde{L}_k}{J}\right)^{u-1}\right]\widetilde{Q}_{(q-u),k} \tag{4.58}$$
$$- q2^{-q}\widetilde{Q}_{(q-1),k} + \widetilde{Q}_{q,k}$$

on $\{\widetilde{L}_k > \alpha\}$.

The last estimate which we need is that

$$\left(\frac{\widetilde{L}_k}{J}\right)^{u-1}\widetilde{Q}_{(q-u),k} \le \widetilde{Q}_{(q-1),k} \qquad \text{for } 2 \le u \le q. \tag{4.59}$$

Before we prove this, we show that it implies the lemma. Indeed, (4.58) and (4.59) show that on $\{\widetilde{L}_k > \alpha\}$,

$$E\{\widetilde{Q}_{q,k+1}|\mathcal{G}_k\} \le \left[D_{12}\sum_{u=2}^{q}\left(\frac{J}{\widetilde{L}_k}\right)^{u-1} + D_{12}q\varepsilon - q2^{-q}\right]\widetilde{Q}_{(q-1),k} + \widetilde{Q}_{q,k}$$
$$\le \left[D_{12}\sum_{u=2}^{q}\left(\frac{J}{\alpha}\right)^{u-1} + D_{12}q\varepsilon - q2^{-q}\right]\widetilde{Q}_{(q-1),k} + \widetilde{Q}_{q,k}.$$

If $\varepsilon$ is chosen so that

$$0 < \varepsilon < \frac{1}{D_{12}2^{q+3}} \tag{4.60}$$

and $\alpha > J$ is chosen so that

$$D_{12}\frac{J}{\alpha - J} \le q2^{-q-2} \quad \text{and} \quad q2^{-q}\frac{\alpha^{q-1}}{J^{q-2}} \ge (D_{13}J) \vee 2,$$

then we have, on $\{\widetilde{L}_k > \alpha\}$,

$$E\{\widetilde{Q}_{q,k+1}|\mathcal{G}_k\} \le -q2^{-q-1}\widetilde{Q}_{(q-1),k} + \widetilde{Q}_{q,k} \le \widetilde{Q}_{q,k} - 1. \tag{4.61}$$

Thus, (4.38) will follow from (4.59). In turn, (4.38) implies that $\nu_1(\alpha)$ is almost surely finite and the chain $\{U_k\}$ is recurrent, by the well-known criterion of Foster [see Theorem 2.2.1 in Fayolle, Malyshev and Menshikov (1995)].

As for (4.39), note that regardless of the value of $\widetilde{L}_k$, we still have by (4.56), (4.59) that

$$E\{\widetilde{Q}_{q,k+1}|\mathcal{G}_k\} - \widetilde{Q}_{q,k} \le C_2\widetilde{Q}_{(q-1),k} + D_{13}J \tag{4.62}$$

for some constant $C_2 = C_2(J,q,\varepsilon)$. If $\widetilde{L}_k \le \alpha$, then $X_j(\tau_k-) \le \alpha + 1$ for each $j$ and hence $\widetilde{Q}_{(q-1),k} \le J(\alpha+1)^{q-1}$. In particular, on the event $\{\nu_1 = \ell\} \in \mathcal{G}_\ell$, we have $\widetilde{L}_\ell \le \alpha$ and

$$E\{\widetilde{Q}_{q,\ell+1}|\mathcal{G}_\ell\} - \widetilde{Q}_{q,\ell} \le C_3 \tag{4.63}$$

for some constant $C_3 = C_3(J, q, \varepsilon, \alpha)$. Now, let $k \geq \ell + 1$, multiply (4.38) by $I[\nu_2(\alpha) > k, \nu_1(\alpha) = \ell]$ and take conditional expectations first with respect to $\mathcal{G}_k$ and then with respect to $\mathcal{G}_\ell$. This gives

$$\begin{aligned} &E\{\widetilde{Q}_{q,k+1} I[\nu_2(\alpha) > k, \nu_1(\alpha) = \ell] | \mathcal{G}_\ell\} \\ &\quad \leq E\{\widetilde{Q}_{q,k} I[\nu_2(\alpha) > k, \nu_1(\alpha) = \ell] | \mathcal{G}_\ell\} \\ &\quad \leq E\{\widetilde{Q}_{q,k} I[\nu_2(\alpha) > k-1, \nu_1(\alpha) = \ell] | \mathcal{G}_\ell\}. \end{aligned}$$

By iteration of this inequality, we obtain

$$\begin{aligned} &E\{\widetilde{Q}_{q,k+1} I[\nu_2(\alpha) > k, \nu_1(\alpha) = \ell] | \mathcal{G}_\ell\} \\ &\quad \leq E\{\widetilde{Q}_{q,\ell+1} I[\nu_1 = \ell] | \mathcal{G}_\ell\} \qquad (\text{since we always have } \nu_2 > \nu_1) \\ &\quad \leq [\widetilde{Q}_{q,\ell} + C_3] I[\nu_1 = \ell] \\ &\quad \leq [J(\alpha+1)^q + C_3] I[\nu_1 = \ell] \qquad [\text{by (4.63) and the lines before it}]. \end{aligned}$$

By (4.63), the inequality between the extreme members here remains valid for $k = \ell$. Taking expectation in this inequality and summing over $\ell = 1, 2, \dots, k$ then gives

$$E\{\widetilde{Q}_{q,k+1} I[\nu_2 > k \geq \nu_1]\} \leq [J(\alpha+1)^q + C_3] P\{\nu_1 \leq k\}. \tag{4.64}$$

Further, by (4.61) or (4.38), because $\widetilde{L}_k > \alpha$ on $\{\nu_1 > k\}$,

$$\begin{aligned} E\{\widetilde{Q}_{q,k+1} I[\nu_2 > k, \nu_1 > k]\} &= E\{\widetilde{Q}_{q,k+1} I[\nu_1 > k]\} \\ &< E\{\widetilde{Q}_{q,k} I[\nu_1 > k]\} \qquad [\text{by (4.61) or (4.38)}] \\ &\leq E\{\widetilde{Q}_{q,k} I[\nu_1 > k-1]\} \\ &\leq \cdots \leq E\{\widetilde{Q}_{q,1} I[\nu_1 > 0]\} = E\{\widetilde{Q}_q, 1\} < \infty. \end{aligned}$$

The last in equality was proven just before (4.50). Adding this to (4.64) finally gives (4.39) for $k \geq 1$. For $k = 0$ (4.39) again reduces to $E\{\widetilde{Q}_{q,1}\} < \infty$, since $\nu_2 > \nu_1 \geq 1$ by definition.

It remains to prove (4.59). To this end, we recall that $[\widetilde{L}_k/J]^{u-1} \leq \widetilde{Q}_{(u-1),k}/J$ for $u \geq 2$ [see (4.57)]. Thus, it suffices to prove

$$\frac{1}{J}\widetilde{Q}_{(u-1),k} \frac{1}{J}\widetilde{Q}_{(q-u),k} \leq \frac{1}{J}\widetilde{Q}_{(q-1),k} \qquad \text{for } 1 \leq u \leq q. \tag{4.65}$$

But this is a simple case of the Harris–FKG inequality (with respect to the measure which assigns mass $1/J$ to each point $X_j(\tau_k-) - I[j = r(k)]$, $1 \leq j \leq J$). $\square$

Proof of Theorem 4. We take $q = 2$ and fix $\varepsilon$ as in (4.60). We also fix $J$ such that $J \geq J(2, \varepsilon)$ [see Lemma 2 for $J(q, \varepsilon)$], $J \geq J(2)$ and $\alpha \geq \alpha(2, J)$ [see Lemma 3 for $J(2), \alpha(2, J)$]. We abbreviate $I[\nu_1 > k]$ to $I_k$ and similarly write $I_k^{(2)}$ for $I[\nu_2 > k]$. As we already pointed out, it suffices to prove that for any initial state, (4.12) and (4.13) hold. We claim that, in turn, these inequalities will follow from

$$\sum_{k=1}^{\infty} [P\{\nu_2 > k\}]^{1/2} < \infty. \tag{4.66}$$

To see this, first note that $P\{\nu_2 > k\} \leq [P\{\nu_2 > k\}]^{1/2}$. Thus, (4.66) will imply (4.12).

Next, recall that we already showed in the beginning of the proof of Lemma 3 that $\tau_k < \infty$ a.s., so the $U_k$ of (4.11) is well defined for all $k$ [see (4.40)]. We also showed just before (4.50) that for any initial state, $E\{\widetilde{Q}_{2,k}\} < \infty$. From (4.19) with $p = 0$ and $q = 1$, we then see that

$$E\{\delta_k\} \leq \varepsilon E\left\{\frac{\widetilde{L}_k}{J}\right\} = \varepsilon E\left\{\frac{\widetilde{Q}_{1,k}}{J}\right\} \leq \varepsilon E\left\{\left[\frac{\widetilde{Q}_{2,k}}{J}\right]^{1/2}\right\} < \infty.$$

We then also have

$$E\{\tau_k\} = E\left\{\sum_{j=0}^{k-1} \delta_j\right\} < \infty. \tag{4.67}$$

Now, assume that (4.66) has been proven and use the relations

$$\tau_{\nu_2} = \tau_1 + \sum_{k=1}^{\infty} \delta_k I_k^{(2)}$$

and

$$\begin{aligned}
E\{\delta_k I_k^{(2)}\} &\leq \frac{\varepsilon}{J} E\{\widetilde{L}_k I_k^{(2)}\} \qquad \text{[by (4.19) with } q = 1,\ p = 0 \text{ and } \{\nu_2 > k\} \in \mathcal{G}_k] \\
&\leq \frac{\varepsilon}{J} [E\{[\widetilde{L}_k]^2 I_k^{(2)}\} P\{\nu_2 > k\}]^{1/2} \\
&\leq C_4 [E\{\widetilde{Q}_{2,k} I_k^{(2)}\} P\{\nu_2 > k\}]^{1/2} \\
&\qquad\qquad \text{[by the Schwarz' inequality or by (4.35)]} \\
&\leq C_4 [C_1 P\{\nu_2 > k\}]^{1/2} \qquad \text{[by } I_k^{(2)} \leq I_{k-1}^{(2)} \text{ and (4.39)]}.
\end{aligned}$$

Moreover, as we just proved, $E\{\tau_1\} < \infty$, so (4.66) will indeed imply (4.13).

Now, to prove that (4.66) indeed holds, we consider the $\{\mathcal{G}_n\}$-martingale

$$M_n := \sum_{k=1}^{n-1} [\widetilde{Q}_{2,k+1} - \widetilde{Q}_{2,k} - E\{\widetilde{Q}_{2,k+1} - \widetilde{Q}_{2,k} | \mathcal{G}_k\}] I_k^{(2)}, \qquad n \geq 1.$$

We note that on $\{\nu_2 > n\}$, each $I_k^{(2)}, 1 \le k \le n$, equals 1, so by virtue of (4.38),

$$E\{\widetilde{Q}_{2,k+1} - \widetilde{Q}_{2,k} | \mathcal{G}_k\} \le -1, \qquad 1 \le k \le n, k \ne \nu_1.$$

On $\{\nu_1 = k\}$,

$$E\{\widetilde{Q}_{2,k+1} - \widetilde{Q}_{2,k} | \mathcal{G}_k\} \le C_3,$$

by virtue of (4.63). Consequently,

$$\begin{aligned} M_n &\ge \sum_{k=1}^{n-1} [\widetilde{Q}_{2,k+1} - \widetilde{Q}_{2,k}] + n - 2 - C_3 \\ &= \widetilde{Q}_{2,n} - \widetilde{Q}_{2,1} + n - 2 - C_3 \qquad \text{on } \{\nu_2 > n\}. \end{aligned}$$

Since $\widetilde{Q}_{2,n} \ge 0$ by definition, it follows that

$$\begin{aligned} P\{\nu_2 > n\} &\le P\{M_n + \widetilde{Q}_{2,1} \ge n - 2 - C_3\} \\ &\le (n - 2 - C_3)^{-5} 2^4 \{E\{|M_n|^5\} + E\{[\widetilde{Q}_{2,1}]^5\}\}. \end{aligned}$$

But, (4.35) and our remarks just before (4.50) show that

$$E\{[\widetilde{Q}_{2,1}]^5\} \le J^4 E\{\widetilde{Q}_{10,1}\} < \infty$$

for large enough $J$. Thus, it suffices for (4.66) to prove

$$E\{|M_n|^5\} < C_5 n^{5/2}. \tag{4.68}$$

But, by Burkholder's inequality [see Gut (1988), Theorem A.2.2 or Hall and Heyde (1980), Theorem 2.10],

$$\begin{aligned} E\{|M_n|^5\} &\le C_6 E\left\{ \left| \sum_{k=1}^{n-1} [\widetilde{Q}_{2,k+1} - \widetilde{Q}_{2,k} - E\{\widetilde{Q}_{2,k+1} - \widetilde{Q}_{2,k} | \mathcal{G}_k\}]^2 I_k^{(2)} \right|^{5/2} \right\} \\ &\le C_7 n^{3/2} E\left\{ \sum_{k=1}^{n-1} |\widetilde{Q}_{2,k+1} - \widetilde{Q}_{2,k} - E\{\widetilde{Q}_{2,k+1} - \widetilde{Q}_{2,k} | \mathcal{G}_k\}|^5 I_k^{(2)} \right\} \\ &\hspace{20em} [\text{by } (4.35)] \\ &\le C_7 3^4 n^{3/2} \sum_{k=1}^{n-1} E\{|\widetilde{Q}_{2,k+1}|^5 I_k^{(2)}\} + C_7 3^4 n^{3/2} \sum_{k=1}^{n-1} E\{|\widetilde{Q}_{2,k}|^5 I_k^{(2)}\} \\ &\quad + C_7 3^4 n^{3/2} \sum_{k=1}^{n-1} E\{|E\{\widetilde{Q}_{2,k+1} - \widetilde{Q}_{2,k} | \mathcal{G}_k\}|^5 I_k^{(2)}\} \qquad [\text{by } (4.35)] \\ &\le C_7 3^4 (1 + 2^4) n^{3/2} \sum_{k=1}^{n-1} E\{|\widetilde{Q}_{2,k+1}|^5 I_k^{(2)}\} \end{aligned}$$

$$+ C_7 3^4(1+2^4)n^{3/2}\sum_{k=1}^{n-1} E\{|\widetilde{Q}_{2,k}|^5 I_k^{(2)}\}$$

$$\le C_8 n^{3/2}\sum_{k=1}^{n} E\{|\widetilde{Q}_{2,k}|^5 I_{k-1}^{(2)}\}.$$

Finally, by (4.35) once more,

$$|\widetilde{Q}_{2,k}|^5 \le J^4 \widetilde{Q}_{10,k}$$

and $E\{\widetilde{Q}_{10,k} I_{k-1}^{(2)}\}$ is bounded in $k$ by virtue of (4.39). Thus, (4.68) and (4.66) hold. This proves (4.13) and, as we pointed out before, it also proves (4.17) and (4.8).

Since (4.8) is our main conclusion in Theorem 4, we leave it to the interested reader to prove that $(X_1(t), \dots, X_J(t))$ is irreducible and positive recurrent. □

**Acknowledgment.** We thank Ronald Dickman for carrying out some simulations of the DLA model for us and for providing us with Figures 1–3.

## REFERENCES

Alves, O. S. M., Machado, F. P. and Popov, S. Yu. (2002). The shape theorem for the frog model. *Ann. Appl. Probab.* **12** 533–546. MR1910638

Chayes, L. and Swindle, G. (1996). Hydrodynamic limits for one-dimensional particle systems with moving boundaries. *Ann. Probab.* **24** 559–598. MR1404521

Chow, Y. S. and Teicher, H. (1978). *Probability Theory Independence Interchangeability Martingales*, 3rd ed. Springer, New York. MR0513230

Chung, K. L. (1967). *Markov Chains with Stationary Transition Probabilities*, 2nd ed. Springer, New York. MR0217872

Doob, J. L. (1953). *Stochastic Processes*. Wiley, New York. MR0058896

Fayolle, G., Malyshev, V. A. and Menshikov, M. V. (1995). *Topics in the Constructive Theory of Countable Markov Chains*. Cambridge Univ. Press. MR1331145

Gut, A. (1988). *Stopped Random Walks*. Springer, New York. MR0916870

Hall, P. and Heyde, C. C. (1980). *Martingale Limit Theory and Its Application*. Academic Press, New York. MR0624435

Kesten, H. (1987). Hitting probabilities of random walks on $\mathbb{Z}^d$. *Stochastic Process. Appl.* **25** 165–184. MR0915132

Kesten, H. and Sidoravicius, V. (2003a). Branching random walk with catalysts. *Electron. J. Probab.* **8** paper #5. MR1961167

Kesten, H. and Sidoravicius, V. (2003b). The spread of a rumor or infection in a moving population. arXiv math.PR/0312496. MR2184100

Kesten, H. and Sidoravicius, V. (2005). The spread of a rumor or infection in a moving population. *Ann. Probab.* **33** 2402–2462. MR2184100

Kesten, H. and Sidoravicius, V. (2006). A phase transition in a model for the spread of an infection. *Illinois J. Math.* **50** 547–634. MR2247840

Lawler, G. F., Bramson, M. and Griffeath, D. (1992). Internal diffusion limited aggregation. *Ann. Probab.* **20** 2117–2140. MR1188055

Ramirez, A. F. and Sidoravicius, V. (2004). Asymptotic behavior of a stochastic combustion growth process. *J. European Math. Soc.* **6** 293–334. MR2060478

Spitzer, F. (1976). *Principles of Random Walk*, 2nd ed. Springer, New York. MR0388547

Voss, R. F. (1984). Multiparticle fractal aggregation. *J. Stat. Phys.* **36** 861–872.

Witten, T. A. and Sander, L. M. (1981). Diffusion-limited aggregation, a kinetic critical phenomenon. *Phys. Rev. Lett* **47** 1400–1403.

Department of Mathematics
Cornell University
Malott Hall
Ithaca, New York 14853
USA
E-mail: kesten@math.cornell.edu

IMPA
Estr. Dona Castorina 110
Jardim Botânico
CEP 22460-320
Rio de Janeiro, RJ
Brasil
E-mail: vladas@impa.br